\documentclass{amsart} 
\usepackage{graphicx}
\usepackage{amssymb, amsmath, sidecap}
\usepackage{amsfonts}
\usepackage{amssymb}
\usepackage{float}
\usepackage{extarrows}
\usepackage{mathrsfs}
\usepackage{booktabs}
\usepackage{verbatim}
\usepackage{hyperref}
\usepackage[usenames,dvipsnames]{xcolor}
\usepackage{esint}

\renewcommand{\baselinestretch}{1}

\newcommand{\N}{\mathbb N}

\newcommand{\R}{\mathbb R}

\newcommand{\dist}{\text{\rm dist}}

\def\bt{\begin{thm}}
\def\et{\end{thm}}
\def\bl{\begin{lem}}
\def\el{\end{lem}}
\def\bd{\begin{defi}}
\def\ed{\end{defi}}
\def\bc{\begin{cor}}
\def\ec{\end{cor}}
\def\bp{\begin{proof}}
\def\ep{\end{proof}}
\def\br{\begin{rem}}
\def\er{\end{rem}}

\def\d{\, \mathrm{d}}

\def\be{\begin{equation}}
\def\ee{\end{equation}}
\def\bes{\begin{equation*}}
\def\ees{\end{equation*}}
\def\bea{\begin{equation} \begin{aligned}}
\def\eea{\end{aligned} \end{equation}}
\def\beas{\begin{equation*} \begin{aligned}}
\def\eeas{\end{aligned} \end{equation*}}
\def\ba{\begin{align}}
\def\ea{\end{align}}
\def\bas{\begin{align*}}
\def\eas{\end{align*}}

\newtheorem{thm}{Theorem}[section]
\newtheorem{lem}{Lemma}[section]
\newtheorem{defi}{Definition}[section]

\newtheorem{prop}[thm]{Proposition}

\newtheorem{rem}{Remark}[section]
\newtheorem{cor}{Corollary}[section]

\numberwithin{equation}{section}
\numberwithin{figure}{section}

\begin{document}
\title{On the Generalized Harmonic Measure\footnotemark[1]}
\author[Liu]{Duchao Liu}
\address[Duchao Liu]{School of Mathematics and Statistics, Lanzhou
University, Lanzhou 730000, P. R. China} \email{liudch@lzu.edu.cn; liuduchao@gmail.com,
Tel.: +8613893289235, fax: +8609318912481}

\author[Wang]{Yunjie Wang}
\address[Yunjie Wang]{School of Mathematics and Statistics, Lanzhou
University, Lanzhou 730000, P. R. China} \email{220220934191@lzu.edu.cn}

\author[Li]{Qiuli Li}
\address[Qiuli Li]{School of Mathematics and Statistics, Lanzhou
University, Lanzhou 730000, P. R. China} \email{qlli@lzu.edu.cn}

\footnotetext[1]{
Research supported by the National Natural Science Foundation of
China (NSFC 12371096).}


\makeatletter
\@namedef{subjclassname@2020}{\textup{2020} Mathematics Subject Classification}
\makeatother

\begin{abstract}
The axiomatization of harmonic measure theory is established, including the generalized maximum principle, Harnack inequality and Harnack principle. As the applications of the established theory, Dahlberg's theory is generalized. The theory provides an interpretation from the measure theory view point of the elliptic equation theory.
\end{abstract}
 \maketitle

\tableofcontents

\section{Introduction}\label{Sec1}

\vspace{0.2cm}

\subsection{The classical harmonic measure}\label{sec1_1}
It is well known that harmonic measure can be defined in several analytic ways, one of which has relationship with generalized Dirichlet problem, see \cite{Keni05, Helm69}. The generalized Dirichlet problem is to find a harmonic function $h$ on $D$ corresponding to the boundary function $f$, where $D$ is any nonempty open subset of $\mathbb{R}^n$ and $f$ is any extended real-valued function on $\partial D$.

\vspace{0.2cm}

In order to give the classical concept of harmonic measure, we need some related conceptions. An extended real-valued function $u$ is called \emph{hyperharmonic} (\emph{hypoharmonic}) on $D$ if it is superharmonic (subharmonic) or identically $+\infty$ ($-\infty$) on each component of $D$. Denote 
\begin{equation*}
\begin{aligned}
\mathcal{U}_f:=\{u:\,u&\text{ is hyperharmonic and bounded below on }D,\\
&\mathop{\lim\inf}\limits_{y\rightarrow x}u(y)\geq f(x), \forall x\in\partial D\},\\
\mathcal{L}_f:=\{u:\,u&\text{ is hypoharmonic and bounded above on }D, \\
&\mathop{\lim\sup}\limits_{y\rightarrow x}u(y)\leq f(x), \forall x\in\partial D\},
\end{aligned}
\end{equation*}
and $\overline H_f:=\inf\{u:u\in\mathcal{U}_f\}$, $\underline{H}_f:=\sup\{u:u\in\mathcal{L}_f\}$. If $\overline H_f=\underline{H}_f$ and both harmonic on $D$, then $f$ is called  \emph{resolutive} and $H_f:=\overline H_f=\underline{H}_f$ is called the \emph{Dirichlet solution} for $f$. The following theorem confirms that the continuous function is resolutive.

\vspace{0.2cm}

\begin{thm}[Wiener]\label{Wiener}
If $f$ is a continuous real-valued function on the boundary $\partial D$ of the bounded open set $D$, then $f$ is resolutive.
\end{thm}

\vspace{0.2cm}

Suppose $D$ is a bounded open set. Denote $L_x(f):=H_f(x)$ for any $f\in C(\partial D)$ and any $x\in D$. Then for each $x\in D$, $L_x(f)$ is a positive functional on $C(\partial D)$ and moreover, by the Riesz representation theorem, there is a unique unit measure $\omega_x$ on the Borel subsets of $\partial D$ such that $H_f(x)=L_x(f)=\int f\,\mathrm{d}\omega_x$. The measure $\omega_x$ on $C(\partial D)$ is called the \emph{harmonic measure} relative to $D$ and $x$.

\vspace{0.2cm}

\subsection{Our aim}
\vspace{0.2cm}
Our aim in the current paper is to give a self-contained harmonic measure theory in the framework of measure theory in Hausdorff spaces. Some motivations can be found in Section \ref{Sec3}. In this theory, the maximum principle holds in open set depending on the regular closeness of the harmonic measure system. Harnack principle holds depending on the translation and scaling invariant of the harmonic measure system and the boundedness of the relative kernel. As the application of maximum principle and Harnack principle, a generalization of Dahlberg's theory is considered. 

\vspace{0.2cm}

Several examples, see in Section \ref{Sec3}, can be viewed as the special case of our theory, including:
\begin{enumerate}
\item[(1)] The harmonic measure system generated by the classical Laplace operator in $\mathbb{R}^n$;
\item[(2)] The harmonic measure system generated by the degenerate elliptic operator of the form $L=-\text{div}A\nabla$ in $\mathbb{R}^n$, where $A$ satisfies ellipticity and boundedness condition with some different homogeneity. In this case, the unbounded and Co-dimension higher than $1$ boundary are included in the harmonic measure system;
\item[(3)] The harmonic measure system generated by the Laplace operator on a graph.
\end{enumerate}

\vspace{0.2cm}

The paper is organized as follows. In Section \ref{Sec1}, we give the classical definition of the harmonic measure in an analytic way. In Section \ref{Sec2}, the axiomatization of harmonic measure is presented, including the maximal principle, Harnack inequality, Harnack principle and some related topics. In Section \ref{Sec3}, some example of our generalized harmonic measure are presented, including: 1) The classic harmonic measure system; 2) A case of unbounded and Co-dimension higher than $1$ boundary harmonic measure system; 3) A generalization of the Dirichlet boundary problem in the graph theory. In Section \ref{Sec4}, as an application of the maximal principle, Harnack inequality and Harnack principle, a generalization of Dahlberg's theory is considered.

\vspace{0.2cm}

\section{The generalized harmonic measure}\label{Sec2}

\vspace{0.2cm}

In this section, we give the basic axiomatization of the harmonic measure. In Subsection \ref{subsec2_1} and Subsection \ref{subsec2_2}, we will present some basic concepts in topological space theory and measure theory in order to introduce the notations we will take in this paper. Our main discussion starts from Subsection \ref{subsec2_3}.

\subsection{Topological space theory}\label{subsec2_1}

We list some basic concepts and conclusions of topological space theory.

\vspace{0.2cm}

Suppose $X$ is a nonempty set. The pair $(X; \tau)$, or briefly denoted by $X$,  is called a \emph{topological space} with the open sets family $\tau\subseteq2^X$ provided
\begin{enumerate}
\item[(1)] $X,\emptyset\in\tau$;
\item[(2)] $\mathop{\cup}\limits_{\lambda\in\Lambda}O_\lambda\in\tau$ for any $\{O_\lambda\}_{\lambda\in\Lambda}\subseteq\tau$;
\item[(3)] $O_\alpha\cap O_\beta\in\tau$ for any $O_\alpha, O_\beta\in\tau$.
\end{enumerate}
Elements of $\tau$ are called the \emph{open sets} of $X$. Any open set containing $x\in X$ is called a \emph{neighborhood} of $x$ and is denoted by $O_x$. Let $\{x_n\}_{n\in\mathbb{N}}\subseteq X$. By $x_n\rightarrow x_0$ we mean for any open set $U\ni x_0$, there exists an $N\in\mathbb{N}$, for any $n>N$, such that $x_n\in U$. Please note that in some topological spaces, for $x_0\not=y_0$, both $x_n\rightarrow x_0$ and $x_n\rightarrow y_0$ may happen at the same time. We say $x\in X$ is a \emph{boundary point} of the subset $U\subseteq X$, if any neighborhood of $x$ both contains both a point of $U$ and a point of $U^c:=X\backslash U$. All boundary points of $U$ is denoted by $\partial U$. A subset $U\subseteq X$ is always considered as a topological space $(U;\tau_U)$ with the open sets family $\tau_U:=\{U\cap O:O\in\tau\}$. We say that $A$ is a \emph{compact set} in $X$ if any open cover of $A$ has a finite subcover.

\vspace{0.2cm}

Let $X$ and $Y$ be two topological spaces. A map $f:X\rightarrow Y$ is called \emph{continuous} if for each open set $U\subseteq Y$, the set $f^{-1}(U)$ is open.

\vspace{0.2cm}

We say a topological space is a \emph{Hausdorff space} if for any $x,y\in X\,(x\not=y)$, there exist separately two neighborhood of $x,y$, denoted by $O_x,O_y\in\tau$, such that $O_x\cap O_y=\emptyset$. In a Hausdorff space, if $x_n\rightarrow x_0$, then $x_0$ must be unique.

\subsection{Measure theory}\label{subsec2_2}

We list some basic concepts and conclusions of measure theory, see \cite{Evan91, Royd88}.

\vspace{0.2cm}

Suppose $X$ is a nonempty set. A mapping $\omega:2^X\rightarrow[0,\infty]$ is called a \emph{measure} on $X$ provided
\begin{enumerate}
\item[(1)] $\omega(\emptyset)=0$;
\item[(2)] if $A\subseteq\mathop{\cup}\limits_{k=1}^{\infty}A_k$, then $\omega(A)\leq\sum\limits_{k=1}^\infty\omega(A_k)$.
\end{enumerate}

\vspace{0.2cm}

For any $A\subseteq X$ the mapping $\delta:2^X\rightarrow\{0,1\}$ defined by
\begin{equation*}
\delta_x(A)=\left\{\begin{aligned}
                      1,\text{ if }x\in A;\\
                      0,\text{ if }x\not\in A,
                       \end{aligned}\right.
\end{equation*}
is a measure on $X$. We always call this measure \emph{Dirac measure}.

\vspace{0.2cm}

A set $A\subseteq X$ is \emph{$\omega$-measurable} if for each set $B\subseteq X$ such that
\begin{equation*}
\omega(B)=\omega(B\cap A)+\omega(B\backslash A).
\end{equation*}
If $\{A_k\}_{k=1}^\infty$ is a sequence of disjoint $\omega$-measurable sets, then
\begin{equation*}
\omega\left(\mathop{\cup}\limits_{k=1}^\infty A_k\right)=\sum_{k=1}^\infty\omega(A_k).
\end{equation*}

\vspace{0.2cm}

A measure $\omega$ on $X$ is called \emph{Borel measure} if each Borel set in $X$ is $\omega$-measurable.

\vspace{0.2cm}

Assume $\mu$ and $\nu$ are measures on $X$. The measure $\nu$ is called \emph{absolutely continuous} with respect to $\mu$, written $\nu\ll\mu$, provided for any $A\subseteq X$, $\mu(A)=0$ implies $\nu(A)=0$; The measure $\nu$ and $\mu$ are called \emph{mutually singular}, written $\nu\perp\mu$, if there exists a $B\subseteq X$ such that $\mu(X\backslash B)=\nu(B)=0$.

\vspace{0.2cm}

Let $X$ be a set and $Y$ be a topological space. A map $f:X\rightarrow Y$ is called \emph{$\omega$-measurable} if for each open set $U\subseteq Y$, the set $f^{-1}(U)$ is $\omega$-measurable. A function $f:X\rightarrow[-\infty,+\infty]$ is $\omega$-measurable if and only if $f^{-1}([-\infty,a))$ is $\omega$-measurable for any $a\in\mathbb{R}$. 

\vspace{0.2cm}

A function $g:X\rightarrow[-\infty,+\infty]$ is called a \emph{simple} function if the image of $g$ is finite. Denote $g^{\pm}:=\max\{\pm g,0\}$. For simple, $\omega$-measurable function $g$, if either $\sum\limits_{0\leq y\leq\infty}y\omega((g^+)^{-1}(y))<\infty$ or $\sum\limits_{0\leq y\leq\infty}y\omega((g^-)^{-1}(y))<\infty$, we define its \emph{integral}
\begin{equation*}
\begin{aligned}
\int g\,\mathrm{d}\omega:&=\int g^+\,\mathrm{d}\omega-\int g^-\,\mathrm{d}\omega\\
&:=\sum_{0\leq y\leq\infty}y\omega((g^+)^{-1}(y))-\sum_{0\leq y\leq\infty}y\omega((g^-)^{-1}(y)),
\end{aligned}
\end{equation*}
where $g^\pm=\max\{\pm g,0\}$. The \emph{integral} of a $\omega$-measurable function $f:X\rightarrow[-\infty,+\infty]$, denoted by $\int f\,\mathrm{d}\omega$, is defined by the value (maybe $\pm\infty$)
\begin{equation*}
\begin{aligned}
&\inf\left\{\int g\,\mathrm{d}\omega:g\text{ is simple, }\omega\text{-measurable},g\geq f,\int g\,\mathrm{d}\omega\text{ exists}\right\}\\
=:&\int^*f\,\mathrm{d}\omega=\int_{*}f\,\mathrm{d}\omega\\
:=&\sup\left\{\int g\,\mathrm{d}\omega:g\text{ is simple, }\omega\text{-measurable},g\leq f,\int g\,\mathrm{d}\omega\text{ exists}\right\}.
\end{aligned}
\end{equation*}
It is well known that the integral of a function $f:X\rightarrow[-\infty,+\infty]$ exists if and only if $f$ is $\omega$-measurable. A function $f:X\rightarrow[-\infty,+\infty]$ is called \emph{$\omega$-integrable} if $f$ is $\omega$-measurable and $\int|f|\,\mathrm{d}\omega<\infty$. For $\delta_x$-measurable function $f$, $\int f\,\mathrm{d}\delta_x=f(x)$.

\vspace{0.2cm}

For $\omega$-measurable function we always denote $\omega(f):=\int f\,\mathrm{d}\omega$, and for measurable set $A\subseteq X$, briefly denote $\omega(A):=\omega(\chi_A)=\int \chi_A\,\mathrm{d}\omega$.

\vspace{0.2cm}

For $1\leq p<\infty$, $L^p(X;\omega)$ (resp. $L_{\text{loc}}^p(X;\omega)$) denotes the set of all $\omega$-measurable functions on $X$ such that $|f|^p$ is integrable (resp. locally integrable). $L^\infty(X;\omega)$ (resp. $L_{\text{loc}}^\infty(X;\omega)$) denotes the set of all $\omega$-measurable functions on $X$ such that $|f|$ is essentially bounded (resp. locally and essentially bounded).

\vspace{0.2cm}

Let $\omega,\omega_k\,(k\in\mathbb{N})$ be measures on a topological space $X$. We say the measures $\{\omega_k\}_{k\in\mathbb{N}}$ \emph{weakly converge} to the measure $\omega$, written $\omega_k\rightharpoonup\omega$ if for all $f\in C_c(X)$,
\begin{equation*}
\lim_{k\rightarrow\infty}\omega_k(f)=\omega(f),
\end{equation*}
where $C_c(X)$ consists of continuous functions on $X$ with compact support.

\vspace{0.2cm}

A measure $\omega$ on a topological space $X$ is called \emph{$\sigma$-finite} if there exists a sequence of $\omega$-measurable sets $\{X_n\}$ such that $X=\mathop{\cup}\limits_{n=1}^\infty X_n$ with $\omega(X_n)<\infty$.

\vspace{0.2cm}

\begin{thm}[Radon-Nikodym Theorem]
Let $\mu$ be a $\sigma$-finite Borel measures on a topological space $X$ and $\nu$ be a Borel measure on $X$ with $\nu\ll\mu$. Then there is a nonnegative $\mu$-measurable function $f$ such that for any Borel set $E\subseteq X$
\begin{equation*}
\nu(E)=\int_Ef\,\d\mu.
\end{equation*}
The function $f$ is unique in the sense of $\mu$-a.e..
\end{thm}

The function $f$ given by Radon-Nikodym Theorem is called the \emph{Radon-Nikodym derivative} of $\nu$ with respect to $\mu$, which is denoted by $\frac{\d\nu}{\d\mu}$.

\vspace{0.2cm}

A measure $\omega$ on $\mathbb{R}^n$ is called \emph{Borel regular} if $\omega$ is Borel measure and for any $A\subseteq\mathbb{R}^n$ there exists a Borel set $B$ such that $A\subseteq B$ and $\omega(A)=\omega(B)$. A measure $\omega$ on $\mathbb{R}^n$ is called \emph{Radon measure} if $\omega$ is Borel regular and $\omega(K)<\infty$ for any compact set $K\subset\mathbb{R}^n$. In this case Radon-Nikodym derivative in Radon-Nikodym Theorem satisfies
\begin{equation*}
\frac{\d\mu}{\d\nu}(x)=\lim_{r\rightarrow0+}\frac{\mu(B(x,r))}{\nu(B(x,r))}.
\end{equation*}

\vspace{0.2cm}

\begin{thm}[Weak compactness for measures in $\mathbb{R}^n$] 
Let $\{\mu_k\}$ be a sequence of Radon measures on $\mathbb{R}^n$ satisfying $\mathop{\sup}\limits_k\mu_k(K)<\infty$ for each compact set $K\subset\mathbb{R}^n$. Then there exists a subsequence $\{\mu_{k_j}\}$ and a Radon measure $\mu$ such that $\mu_{k_j}\rightharpoonup\mu$.
\end{thm}

\vspace{0.2cm}

\subsection{The generalized harmonic measure and maximum principle}\label{subsec2_3}

\vspace{0.2cm}

We give the definition of the generalized harmonic measure defined on the Hausdorff space. 

\vspace{0.2cm}

In this article, we always suppose $X$ is a Hausdorff space.

\vspace{0.2cm}

\begin{defi}[Pre-harmonic measure]\label{pre_HM}
Suppose $D\subset X$ is open with $\partial D\not=\emptyset$. A family of Borel measures $\{\omega^D_x:2^{\partial D}\rightarrow[0,1]:\,x\in D\}$ on $\partial D$, briefly denoted by $\{\omega^D_x\}$, is called a \emph{pre-harmonic measure} relative to $D$ provided for any connected component $V$ of $D$,
\begin{enumerate}
\item[(1)] For any $x\in V$, $\omega^D_x(\partial V)=1$;
\item[(2)] For any $x,y\in V$, $\omega^D_x\ll\omega^D_y$.
\end{enumerate}
\end{defi}

\vspace{0.2cm}

From (2) we get the following proposition.

\vspace{0.2cm}

\begin{prop}
Suppose $\{\omega^D_x\}$ is a pre-harmonic measure relative to $D$. Then the class of $\omega^D_x$-measure $0$ Borel subsets of $\partial D$ is independent of $x\in D$.
\end{prop}

\vspace{0.2cm}

For different $x$ and $y$ in $D$, the measurable sets and functions of $\omega_x^D$ and $\omega_y^D$ may be separately different. We want to find common measurable sets and functions of $\{\omega_x^D:x\in D\}$. Denote 
\begin{equation*}
\mathcal{F}_x:=\{E\Delta N: E\subseteq X \text{ is Borel, } N\subseteq X \text{ is }\omega^D_x\text{-measure }0\}.
\end{equation*}
Then $\mathcal{F}_x$ is a $\sigma$-algebra. Let $\mathcal{F}:=\mathop{\cap}\limits_{x\in D}\mathcal{F}_x$. Then $\mathcal{F}$ is a $\sigma$-Algebra containing Borel subsets of $\partial D$. Since $\omega^D_x$ can be uniquely extended to $\mathcal{F}_x$, the same is true of $\mathcal{F}$. The extension of $\omega^D_x$ to $\mathcal{F}$ is also denoted by $\omega^D_x$.

\begin{defi}[Harmonic measure]\label{H_m}
(1) Elements in $\mathcal{F}$ are called \emph{$\omega_\cdot^D$-measurable sets} and extended real-valued functions measurable relative to $\mathcal{F}$ are called \emph{$\omega_\cdot^D$-measurable functions}. The pre-harmonic measure $\{\omega_x^D\}$ restricted on $\mathcal{F}$, briefly denoted by $\omega_\cdot^D$, is called the \emph{harmonic measure} relative to $D$. 

(2) A $\omega_\cdot^D$-measurable function $f$ on $\partial D$ is called \emph{$\omega_\cdot^D$-integrable} if it is integrable relative to $\omega_x^D$ for any $x\in D$. 

(3) A harmonic measure $\omega_\cdot^D$ is called \emph{continuous} in $D$ if $h_f(x):=\omega_x^D(f)$ is continuous in $D$ for any $f\in L(\partial D;\omega_\cdot^D)$.

(4) A point $y\in\partial D$ is called a \emph{regular point} relative to $\omega_\cdot^D$ if for any $f\in L(\partial D;\omega_\cdot^D)$ continuous at $y\in\partial D$,
\begin{equation}\label{s1}
\lim_{D\ni x_n\rightarrow y}\omega_{x_n}^D(f)=f(y).
\end{equation} 
$\omega_\cdot^D$ or $D$ is called \emph{regular} if every point of $\partial D$ is regular relative to $\omega_\cdot^D$.
\end{defi}

\vspace{0.2cm}

For any $f\in L(\partial D;\omega_\cdot^D)$ and any $y\in\partial D$, we always denote $\omega_{y}^D(f):=f(y)$.

\vspace{0.2cm}

We list several simple properties in the following proposition.

\vspace{0.2cm}

\begin{prop}\label{R_point}
(1) $C(\partial D)\subseteq L(\partial D;\omega_\cdot^D)$;

(2) Suppose $V$ is a component of $D$.  Then the class of $\omega^D_x$-integrable functions is independent of $x\in V$;

(3) Suppose $y\in\partial D$ is a regular point. Then $\omega_{x_n}^D\rightharpoonup\delta_y$ as $D\ni x_n\rightarrow y$, that is, for any $f\in C_c(\partial D)$,
\begin{equation}\label{w1}
\lim_{D\ni x_n\rightarrow y}\omega_{x_n}^D(f)=f(y);
\end{equation}

(4) Suppose $y\in\partial D$ is a regular point. Then for any open set $V\ni y$ of $\partial D$, $\lim\limits_{D\ni x_n\rightarrow y}\omega_{x_n}^D(V)=1$.
\end{prop}

\vspace{0.2cm}

\begin{rem}
We only know the harmonic measure induced by the Laplace operator in $\mathbb{R}^n$ satisfies the equivalence of \eqref{s1} and \eqref{w1}, see Lemma \ref{lem35} or Corollary 8.21 in \cite{Helm69}. It is interesting to derive the condition to insure that \eqref{s1} and \eqref{w1} are equivalent in general topological spaces.
\end{rem}

\vspace{0.2cm}

\begin{prop}\label{H_A}
Suppose $\omega_\cdot^D$ is a harmonic measure. Then, for any nonempty open set $A\subseteq\partial D$ with a regular point in $A$, $\omega_\cdot^D(A)\not=0$.
\end{prop}

\vspace{0.2cm}

\begin{proof}
If there exists a nonempty open set $A$ with a regular point $y\in A$ such that $\omega_\cdot^D(A)=0$, then $\omega_x(\chi_A)=\omega_x(A)=0$ for any $x\in D$. Since $\omega_x(\chi_A)=\omega_x(A)=0$ is continuous with respect to $x\in D$, we get
\begin{equation*}
\lim_{D\ni x_n\rightarrow y}\omega_{x_n}(A)=0.
\end{equation*} 
At the same time, by Definition \ref{H_m} (4), since $y\in A$ is regular and $\chi_A$ is continuous at $y\in A$, we get
\begin{equation*}
\lim_{D\ni x_n\rightarrow y}\omega_{x_n}(A)=\chi_A(y)=1,
\end{equation*}
which contradicts to the above equation.
\end{proof}

\vspace{0.2cm}

Now we consider several different forms maximum principle of the harmonic measure.

\vspace{0.2cm}

\begin{thm} [Strong maximum principle for $\chi_A$ on regular open set]
Suppose $D$ is regular and $A\subseteq\partial D$ is open. Then
\begin{equation*}
\sup_{x\in D}\omega_x^D(\chi_A)=\mathop{\max}\limits_{y\in\partial D}\chi_A(y),
\end{equation*}
and if there exists an $x_0\in D$ such that
\begin{equation*}
\omega_{x_0}(\chi_A)=\sup_{x\in D}\omega_x^D(\chi_A),
\end{equation*}
then $\omega_{x}(\chi_A)$ is constant for $x\in D$.
\end{thm}

\vspace{0.2cm}

\begin{proof}
When $A = \emptyset$, the conclusion holds clearly. 

When $A\subseteq\partial D$ is a nonempty open set, by Proposition \ref{H_A}, we get
\begin{equation*}
0<\omega_x^D(\chi_A)=\omega_x^D(A)\leq1.
\end{equation*}
And for any  $y\in A$, $y$ is regular and $\chi_A$ is continuous at $y\in A$, we get
\begin{equation*}
\lim_{D\ni x_n\rightarrow y}\omega_{x_n}(A)=\chi_A(y)=1.
\end{equation*}
Then
\begin{equation*}
\sup_{x\in D}\omega_x^D(\chi_A)=\mathop{\max}\limits_{y\in\partial D}\chi_A(y)=1.
\end{equation*}
If there exists an $x_0\in D$ such that
\begin{equation*}
\omega_{x_0}(\chi_A)=\sup_{x\in D}\omega_x^D(\chi_A) = 1,
\end{equation*}
then $\omega_{x_0}(\partial D\backslash A)=0$. By Definition \ref{pre_HM} (2) for any $x\in D$, $\omega_{x}(\partial D\backslash A)=0$. This implies that for any $x\in D$, $\omega_{x}(\chi_A)=\omega_{x}(A)=1$.
\end{proof}

\vspace{0.2cm}

\begin{thm}[Weak maximum principle for regular harmonic measure]\label{WMPr}
Suppose $\omega_\cdot^D$ is regular. If $f\in C(\partial D)$ and there exists a $y_0\in\partial D$ such that $f(y_0)=\mathop{\sup}\limits_{y\in\partial D}f(y)$, then
\begin{equation*}
\sup_{x\in D}\omega_x^D(f)=\mathop{\max}\limits_{y\in\partial D}f(y).
\end{equation*}
\end{thm}

\vspace{0.2cm}

\begin{proof}
By (1) in Definition \ref{pre_HM}, for all $x\in D$, it is clear the flowing equation holds
\begin{equation*}
\omega_x^D(f)=\int_{\partial D} f\,\mathrm{d}\omega_x^D\leq\mathop{\max}\limits_{y\in\partial D}f(y)\int_{\partial D}\,\mathrm{d}\omega_x^D=\mathop{\max}\limits_{y\in\partial D}f(y).
\end{equation*}
Then we get
\begin{equation*}
\sup_{x\in D}\omega_x^D(f)\leq\mathop{\max}\limits_{y\in\partial D}f(y).
\end{equation*}
Since $f(y_0)=\mathop{\max}\limits_{y\in\partial D}f(y)$, by Proposition \ref{R_point} (3), $\lim\limits_{x_n\rightarrow y_0}\omega_{x_n}^D(f)=f(y_0)$. Then
\begin{equation*}
\sup_{x\in D}\omega_x^D(f)=\mathop{\max}\limits_{y\in\partial D}f(y).
\end{equation*}
\end{proof}

\vspace{0.2cm}

\begin{thm}[Boundary strong maximum principle]\label{BSMP}
Suppose $f\in L(\partial D;\omega_\cdot^D)\cap L^\infty(\partial D;\omega_\cdot^D)$. If there exists an $x_0\in D$ such that
\begin{equation*}
\omega_{x_0}(f)=\mathop{\mathrm{ess\,sup}}\limits_{y\in\partial D}f(y),
\end{equation*}
then $\omega_{x}(f)$ is constant for $x\in D$.
\end{thm}

\vspace{0.2cm}

\begin{proof}
Denote $c_0=\mathop{\mathrm{ess\,sup}}\limits_{y\in\partial D}f(y)$, $A_0:=\{y\in\partial D:f(y)=c_0\}$, $A_1:=\{y\in\partial D:f(y)>c_0\}$ and $A_2:=\{y\in\partial D:f(y)<c_0\}$. It is clear that $\partial D=A_0\cup A_1\cup A_2$, $A_i(i=1,2,3)$ are disjoint and $\omega_\cdot^D(A_1)=0$. We will prove $\omega_\cdot^D(A_2)=0$. In fact by $\omega_\cdot^D(A_1)=0$ we get
\begin{equation*}
\begin{aligned}
c_0&=\omega_{x_0}(f)=\int_{\partial D}f\,\mathrm{d}\omega_{x_0}^D=\int_{\partial D\backslash A_1}f\,\mathrm{d}\omega_{x_0}^D\\
&=\int_{\partial D\backslash A_1}[c_0+(f-c_0)]\,\mathrm{d}\omega_{x_0}^D\\
&=c_0+\int_{\partial D\backslash A_1}(f-c_0)\,\mathrm{d}\omega_{x_0}^D.
\end{aligned}
\end{equation*}
Then
\begin{equation*}
\int_{A_2}(f-c_0)\,\mathrm{d}\omega_{x_0}^D=\int_{A_2\cup A_0}(f-c_0)\,\mathrm{d}\omega_{x_0}^D=\int_{\partial D\backslash A_1}(f-c_0)\,\mathrm{d}\omega_{x_0}^D=0.
\end{equation*}
Since $f<c_0$ on $A_2$, we get $\omega_\cdot^D(A_2)=0$. Then we can get 
\begin{equation*}
\omega_\cdot^D(A_0)=\omega_\cdot(\partial D)-\omega_\cdot^D(A_1)-\omega_\cdot^D(A_2)=1,
\end{equation*}
which implies $f(y)=c_0$ for $\omega_\cdot^D$-a.e. $y\in\partial D$. Then for any $x\in D$,
\begin{equation*}
\omega_x(f)=\int_{\partial D} f\,\mathrm{d}\omega_x=\int_{A_0} f\,\mathrm{d}\omega_x=c_0.
\end{equation*} 
\end{proof}

\vspace{0.2cm}

By Theorem \ref{WMPr} and Theorem \ref{BSMP}, we can get the following strong maximum principle.

\vspace{0.2cm}

\begin{thm}[Strong maximum principle for regular harmonic measure]\label{SMPr}
Suppose $\omega_\cdot^D$ is regular. If $f\in C(\partial D)$ and there exists a $y_0\in\partial D$ such that $f(y_0)=\mathop{\sup}\limits_{y\in\partial D}f(y)$, then
\begin{equation*}
\sup_{x\in D}\omega_x^D(f)=\mathop{\max}\limits_{y\in\partial D}f(y),
\end{equation*}
and if there exists an $x_0\in D$ such that
\begin{equation*}
\omega_{x_0}(f)=\sup_{x\in D}\omega_x^D(f),
\end{equation*}
then $\omega_{x}(f)$ is constant for $x\in D$.
\end{thm}

\vspace{0.2cm}

For different open sets in $X$, we need the following concept of harmonic measure system.

\vspace{0.2cm}

\begin{defi}[Harmonic measure system]
(1) We call the family of continuous harmonic measures $\{\omega_\cdot^{D_\lambda}: \lambda\in\Lambda\}$ is a \emph{harmonic measure system} on $X$, briefly denoted by $\{\omega_\cdot^{D_{\lambda}}\}$, if its elements are \emph{compatible}, that is, for any $D, \widetilde D\in\{D_{\lambda}\}$ with $\widetilde{D}\cap D\not=\emptyset$, any $f\in L(\partial D;\omega_\cdot^D)$ and any $\widetilde f\in L(\partial \widetilde D;\omega_\cdot^{\widetilde D})$ with $\omega_y^{\widetilde{D}}(\widetilde f)=f(y)$ for any $y\in\widetilde D\cap\partial D$ and $\omega_y^{D}(f)=\widetilde f(y)$ for any $y\in\partial\widetilde D\cap\bar D$, such that $\omega^{D}_x(f)=\omega_x^{\widetilde{D}}(\widetilde f)$ for any $x\in \widetilde{D}\cap D$.

(2) In a harmonic measure system $\{\omega_\cdot^{D_{\lambda}}\}$, a point $y\in\partial D$ is called \emph{regular accessible (RA)} in $D\in\{D_\lambda\}$ if there exists a regular $\omega_\cdot^{\widetilde D}\in\{\omega_\cdot^{D_{\lambda}}\}$ with $\widetilde D\subseteq D$ such that $y\in\partial\widetilde D\cap\partial D$.

(3) A harmonic measure system $\{\omega_\cdot^{D_{\lambda}}\}$ on $X$ is called \emph{regular closed} if for any $D\in\{D_\lambda\}$, there exists an increasing regular sequence $\{D_k\}\subseteq\{D_\lambda\}$ such that $D=\mathop{\cup}\limits_{k=1}^\infty D_k$.
\end{defi}

\vspace{0.2cm}

In a harmonic measure system, the following propositions holds.

\vspace{0.2cm}

\begin{prop}
Suppose in a harmonic measure system $\{\omega_\cdot^{D_{\lambda}}\}$, $D, \widetilde D\in\{D_{\lambda}\}$ with $\widetilde{D}\subseteq D$. Then for any functions $f\in L(\partial D;\omega_\cdot^D)$ and any $x\in \widetilde{D}$, $\omega^{D}_x(f)=\omega_x^{\widetilde{D}}(\omega_{\cdot}^{D}(f))$.
\end{prop}

\vspace{0.2cm}

\begin{thm}\label{R_Comp}
Suppose the following conditions hold: (1) The harmonic measure system $\{\omega_\cdot^{D_\lambda}\}$ in $\mathbb{R}^n$ is translation closed and invariant, that is, for any $x_0\in\mathbb{R}^n$ and $D\in\{D_\lambda\}$, such that $D+\{x_0\}\in\{D_\lambda\}$; for any $x\in D$ and any $A\subseteq\partial D$, such that $\omega_{x+x_0}^{D+\{x_0\}}(A+\{x_0\})=\omega_{x}^D(A)$; (2) $\{D_\lambda\}$ is closed about the finite operation ``$\cup$''; (3) There exists a $\{\omega_\cdot^{O_k}\}$-regular open set sequence $\{O_k\}\subseteq\{D_\lambda\}$ in $\mathbb{R}^n$ such that $\lim\limits_{k\rightarrow\infty}\mathrm{diam}\,O_k=0$. Then $\{\omega_\cdot^{D_\lambda}\}$ is regular closed.
\end{thm}

\vspace{0.2cm}

\begin{proof}
Suppose $D\in\{D_\lambda\}$. Let $\{K_l\}$ be an increasing sequence of compact sets such that $D=\mathop{\cup}\limits_l K_l$. Then there exists a finite covering $\{O^l_i\}_{i=1}^{k_l}$ of $K_l$, where $O^l_i\in\{D_\lambda\}$, $i=1,2,\cdots,k_l$ are regular and $\bar O^l_i\subset D$. By Proposition \ref{cupR} $D_l:=\mathop{\cup}\limits_{j=1}^{k_l}O^l_j\in\{D_\lambda\}$ is regular. It is clear that $\bar D_l\subset D$ and $D=\mathop{\cup}\limits_{l=1}^\infty D_l$.
\end{proof}

\vspace{0.2cm}

\begin{thm}[Strong maximum principle for open set]\label{Str_MP_G}
Suppose the harmonic measure system $\{\omega_\cdot^{D_\lambda}\}$ is regular closed, $D\in\{D_\lambda\}$ and $f\in L(\partial D,\omega_\cdot^D)\cap L^\infty(\partial D,\omega_\cdot^D)$. If there exists an $x_0\in D$ such that
\begin{equation*}
\omega^D_{x_0}(f)=\sup_{x\in D}\omega_x^D(f),
\end{equation*} 
then $\omega_x^D(f)$ is constant for $x\in D$.
\end{thm}

\vspace{0.2cm}

\begin{proof}
Since the harmonic measure system $\{\omega_\cdot^{D_\lambda}\}$ is regular closed, there exists an increasing regular open set sequence $\{D_k\}$ such that $D=\mathop{\cup}\limits_{k=1}^\infty D_k$. By $x_0\in D$, we can find a $\widetilde K\in\mathbb{N}$, for any $k\geq\widetilde K$ such that $x_0\in D_k$. 

For any $x\in D$, we prove $\omega_{x}^D(f)=\omega_{x_0}^D(f)$. In fact, since $x\in D$, there exists a $K\geq\widetilde K$ such that $x_0, x\in D_K$. Then by the compatibility of the harmonic measure system $\{\omega_\cdot^D\}$, we get
\begin{equation*}
\omega_{x_0}^{D_K}(\omega_\cdot^D(f))=\omega_{x_0}^D(f)=\sup_{x\in D}\omega_x^D(f)\mathop{=}\limits^{x_0\in D_K}\sup_{x\in D_K}\omega_x^D(f)=\sup_{x\in D_K}\omega_x^{D_K}(\omega_\cdot^D(f)).
\end{equation*} 
Since $D_K$ is regular and $\omega_\cdot^D(f)$ is continuous on $\partial D_K$, by Theorem \ref{SMPr}, we get
\begin{equation*}
\omega_{x}^D(f)=\omega_x^{D_K}(\omega_\cdot^D(f))=\omega_{x_0}^{D_K}(\omega_\cdot^D(f))=\omega_{x_0}^D(f),
\end{equation*}
which completes the proof.
\end{proof}

\vspace{0.2cm}

From the proof of the above theorem, we can get the following corollary.

\vspace{0.2cm}

\begin{thm}
Suppose $\{\omega_\cdot^{D_\lambda}\}$ is a harmonic measure system and $D\in \{D_\lambda\}$, there exists an increasing regular open set sequence $\{D_k\}$ such that $D=\mathop{\cup}\limits_{k=1}^\infty D_k$ and $f\in L(\partial D,\omega_\cdot^D)\cap L^\infty(\partial D,\omega_\cdot^D)$. If there exists an $x_0\in D$ such that
\begin{equation*}
\omega^D_{x_0}(f)=\sup_{x\in D}\omega_x^D(f),
\end{equation*} 
then $\omega_x^D(f)$ is constant for $x\in D$.
\end{thm}

\vspace{0.2cm}

In a harmonic measure system, the union of two regular sets is regular.

\vspace{0.2cm}

\begin{prop}\label{cupR}
Suppose in a harmonic measure system $\{\omega_\cdot^{D_{\lambda}}\}$, $D_1,D_2,D_1\cup D_2\in\{D_\lambda\}$. If $D_1$ and $D_2$ are regular, then $D_1\cup D_2$ is regular.
\end{prop}

\vspace{0.2cm}

\begin{proof}
Since $D_1$ and $D_2$ are regular, for any $f_1\in L(\partial D_1;\omega_\cdot^{D_1})$ continuous at $y_1\in\partial D_1$, $f_2\in L(\partial D_2;\omega_\cdot^{D_2})$ continuous at $y_2\in\partial D_2$,
\begin{equation*}
\lim_{D_1\ni x_n\rightarrow y_1}\omega_{x_n}^{D_1}(f_1)=f_1(y_1),\quad \lim_{D_2\ni x_n\rightarrow y_2}\omega_{x_n}^{D_2}(f_2)=f_2(y_2).
\end{equation*}

If $D_1\cap D_2=\emptyset$ the conclusion holds obviously. 

If $D_1\cap D_2\not=\emptyset$, then for any $y\in\partial(D_1\cup D_2)$, $y\in\partial D_1\backslash\partial D_2$ or $y\in\partial D_2\backslash\partial D_1$ or $y\in\partial D_1\cap\partial D_2$. 

If $y\in\partial D_1\backslash\partial D_2$, taking the above $f_1(\cdot):=\omega_\cdot^{D_1\cup D_2}(f)$ on $\partial D_1$ for any $f\in L(\partial(D_1\cup D_2);\omega_\cdot^{D_1\cup D_2})$ continuous at $y\in\partial D_1\backslash\partial D_2$, then $f_1\in L(\partial D_1;\omega_\cdot^{\partial D_1})$ is continuous at $y\in\partial D_1$. Then by the compatibility of $\{\omega_\cdot^D\}$, we get
\begin{equation*}
\lim_{D_1\cup D_2\ni x_n\rightarrow y\in\partial D_1\backslash\partial D_2}\omega_{x_n}^{D_1\cup D_2}(f)=\lim_{D_1\ni x_n\rightarrow y\in\partial D_1}\omega_{x_n}^{D_1}(f_1)=f_1(y)=f(y),
\end{equation*}
which means $y\in\partial D_1\backslash\partial D_2$ is regular relative to $D_1\cup D_2$.

If $y\in\partial D_2\backslash\partial D_1$, similar argument of the above also applies. 

If $y\in\partial D_1\cap\partial D_2$, denoting $\{x_{n_1}\}=\{x_n\}\cap (D_1\backslash\bar D_2)$, $\{x_{n_2}\}=\{x_n\}\cap (D_2\backslash\bar D_1)$ and $\{x_{n_3}\}=\{x_n\}\cap (\bar D_1 \cap \bar D_2)$, then for any $f\in L(\partial(D_1\cup D_2);\omega_\cdot^{D_1\cup D_2})$ continuous at $y\in\partial D_1\cap\partial D_2$, similar argument yields
\begin{equation*}
\begin{aligned}
\lim_{(D_1\backslash D_2)\ni x_{n_1}\rightarrow y\in\partial D_1\cap\partial D_2}\omega_{x_{n_1}}^{D_1\cup D_2}(f)=\lim_{D_1\ni x_{n_1}\rightarrow y\in\partial D_1}\omega_{x_{n_1}}^{D_1}(f_1)=f_1(y)=f(y),\\
\lim_{(D_2\backslash D_1)\ni x_{n_2}\rightarrow y\in\partial D_2\cap\partial D_1}\omega_{x_{n_2}}^{D_1\cup D_2}(f)=\lim_{D_2\ni x_{n_2}\rightarrow y\in\partial D_2}\omega_{x_{n_2}}^{D_2}(f_2)=f_2(y)=f(y),\\
\lim_{(D_1\cap D_2)\ni x_{n_3}\rightarrow y\in\partial D_1\cap\partial D_2}\omega_{x_{n_3}}^{D_1\cup D_2}(f)=\lim_{D_1\ni x_{n_3}\rightarrow y\in\partial D_1}\omega_{x_{n_3}}^{D_1}(f_1)=f_1(y)=f(y).
\end{aligned}
\end{equation*}
The above three equations imply
\begin{equation*}
\lim_{D_1\cup D_2\ni x_{n}\rightarrow y\in\partial D_1\cap\partial D_2}\omega_{x_{n}}^{D_1\cup D_2}(f)=f(y),
\end{equation*}
which means $y\in\partial D_1\cap\partial D_2$ is regular relative to $D_1\cup D_2$.
\end{proof}

\vspace{0.2cm}

\begin{rem}
If $D_1$ and $D_2$ are regular in a harmonic measure system $\{\omega_\cdot^{D_{\lambda}}\}$, we can not get $D_1\cap D_2$ is regular in general provided $D_1\cap D_2\in\{D_\lambda\}$.
\end{rem}

\vspace{0.2cm}

We can construct a new harmonic measure $\omega_\cdot^D$ with $D:=D_1\cup D_2$, where $\omega_\cdot^{D_i}, i=1,2,$ are regular in harmonic measure system $\{\omega_\cdot^{D_\lambda}\}\not\ni\omega_\cdot^D$.

\vspace{0.2cm}

Suppose the $\{\omega_\cdot^{D_\lambda}\}$ is a harmonic measure system and $D_1,D_2\in\{D_\lambda\}$ are regular. For any $A\in\partial(D_1\cup D_2)$, we will define a measure $\omega^{D_1\cup D_2}_x$ for any $x\not\in\partial(D_1\cup D_2)$, which is expected to be a harmonic measure.

\vspace{0.2cm}

Without loss of generality, we suppose $D_1\cap D_2\not=\emptyset$. And for any $A\subseteq\partial(D_1\cup D_2)$, we can make decomposition $A=A_1\cup A_2$ with $A_1\cap A_2=\emptyset$, where $A_1:=A\cap\partial D_1$ and $A_2:=A\backslash A_1$. Define the boundary function on $\partial D_2$ by
\begin{equation*}
f^{(2)}_1(y):=\left\{\begin{aligned}
&\omega_y^{D_1}(A_1), && y\in\partial D_2\cap\bar D_1;\\
&\chi_{A_2}(y), &&y\in\partial D_2\backslash\bar D_1.
\end{aligned}\right.
\end{equation*}
By the strong maximum principle Theorem \ref{SMPr} in $D_1$, $0< f^{(2)}_1(y)<1$ for any $y\in\partial D_2\cap D_1$. Define the boundary function on $\partial D_1$ by
\begin{equation*}
f^{(1)}_1(y):=\left\{\begin{aligned}
&\omega_y^{D_2}(f^{(2)}_1), && y\in\partial D_1\cap D_2;\\
&\chi_{A_1}(y), &&y\in\partial D_1\backslash D_2.
\end{aligned}\right.
\end{equation*}
By the strong maximum principle Theorem \ref{SMPr} in $D_2$, $0< f^{(1)}_1(y)<1$ for any $y\in\partial D_1\cap D_2$. By induction, for $i\in\mathbb{N}$, define the boundary function on $\partial D_2$ by
\begin{equation*}
f^{(2)}_i(y):=\left\{\begin{aligned}
&\omega_y^{D_1}(f^{(1)}_{i-1}), && y\in\partial D_2\cap\bar D_1;\\
&\chi_{A_2}(y), &&y\in\partial D_2\backslash\bar D_1.
\end{aligned}\right.
\end{equation*}
And define the boundary function on $\partial D_1$ by
\begin{equation*}
f^{(1)}_i(y):=\left\{\begin{aligned}
&\omega_y^{D_2}(f^{(2)}_i), && y\in\partial D_1\cap D_2;\\
&\chi_{A_1}(y), &&y\in\partial D_1\backslash D_2.
\end{aligned}\right.
\end{equation*}
By the strong maximum principle Theorem \ref{SMPr} in $D_1$ and $D_2$ respectively, we get $f^{(2)}_{i}(y)< f^{(2)}_{i+1}(y)\leq 1$ for any $y\in\partial D_2$ and $f^{(1)}_{i}(y)< f^{(1)}_{i+1}(y)\leq 1$ for any $y\in\partial D_1$. Then $\{f^{(k)}_i:i\in\mathbb{N}\},\,k=1,2$, are point wise increasing and bounded above sequences. This confirms the convergence of $\{f^{(k)}_i\}$. Denote
\begin{equation*}
f^{(1)}(y):=\lim_{i\rightarrow\infty}f^{(1)}_i(y)
\end{equation*}
for any $y\in\partial D_1$ and
\begin{equation*}
f^{(2)}(y):=\lim_{i\rightarrow\infty}f^{(2)}_i(y)
\end{equation*}
for any $y\in\partial D_2$.
Now we define
\begin{equation*}
\omega_x^{D_1\cup D_2}(A):=\left\{\begin{aligned}
&\omega_x^{D_1}(f^{(1)}), && x\in D_1;\\
&\omega_x^{D_2}(f^{(2)}), && x\in D_2.
\end{aligned}\right.
\end{equation*}
Notice that for any $x\in D_1\cap D_2$, by the compatibility of the harmonic measure system,
\begin{equation*}
\omega_x^{D_1}(f^{(1)})=\omega_x^{D_2}(f^{(2)});
\end{equation*}
$f^{(2)}(y)=\omega_y^{D_1}(f^{(1)})$ is continuous for any $y\in\partial D_2\cap\bar D_1$ and $f^{(1)}(y)=\omega_y^{D_2}(f^{(2)})$ is continuous for any $y\in\partial D_1\cap D_2$. 

\vspace{0.2cm}

\begin{thm}
Suppose in a harmonic measure system $\{\omega_\cdot^{D_{\lambda}}\}$, $D_1,D_2\in\{D_\lambda\}$ are regular and $\omega_\cdot^{D_1\cup D_2}$ is defined as the above. Then $\{\omega_\cdot^{D_1\cup D_2}\}\cup\{\omega_\cdot^{D_\lambda}\}$ is a harmonic measure system. 
\end{thm}

\vspace{0.2cm}

\begin{proof}
In fact, we need to prove by steps:
\begin{enumerate}
\item[(a)] $\omega^{D_1\cup D_2}_\cdot$ is a harmonic measure;
\item[(b)] $\{\omega_\cdot^{D_1\cup D_2}\}\cup\{\omega_\cdot^{D_\lambda}\}$ is a harmonic measure system.
\end{enumerate}
We continue to use the notations as the above.

(a) We will prove $\omega^{D_1\cup D_2}_\cdot$ is a harmonic measure.

(a1) We will prove if $A_1\cup A_2=\partial(D_1\cup D_2)$, then $f^{(1)}=f^{(2)}=1$, which implies $\omega^{D_1\cup D_2}_\cdot(\partial(D_1\cup D_2))=1$. We argue by contradiction, if there exists a $y_1\in\partial D_1\backslash D_2$, such that $f^{(1)}(y_1)<1$, then by the strong maximum principle Theorem \ref{SMPr} there exists a $y_2\in\partial D_2\backslash\bar D_1$ such that $f^{(2)}(y_2)<1$. Without loss of generality, we suppose $f^{(1)}(y_1)=\min\limits_{y\in\partial D_1}f^{(1)}(y)$ and $f^{(2)}(y_2)=\min\limits_{y\in\partial D_2}f^{(2)}(y)$. By Theorem \ref{SMPr} in $D_1$ and $D_2$, we get $f^{(1)}(y_1)=f^{(2)}(y_2)$. By the strong maximum principle Theorem \ref{SMPr} in $D_1$ and $D_2$, $f^{(1)}$ and $f^{(2)}$ are constants $\min\limits_{y\in\partial D_2}f^{(2)}(y)<1$ which contradicts to $f^{(1)}(y)=1$ for $y\in A_1$ or $f^{(2)}(y)=1$ for $y\in A_2$.

(a2) We will prove for any $x_1,x_2\in D_1\cup D_2$, $\omega_{x_1}^{D_1\cup D_2}\ll\omega_{x_2}^{D_1\cup D_2}$. In fact, for any $A\subset\partial D_1\cup D_2$ with $\omega_{x_2}^{D_1\cup D_2}(A)=0$, we argue in 4 cases: (1) If $x_1, x_2\in D_2$, then $\omega_{x_2}^{D_2}(f^{(2)})=\omega_{x_2}^{D_1\cup D_2}(A)=0$. By $\omega_{x_1}^{D_2}\ll\omega_{x_2}^{D_2}$, we get $\omega_{x_1}^{D_1\cup D_2}(A)=\omega_{x_1}^{D_2}(f^{(2)})=0$. (2) If $x_1, x_2\in D_1$, this case is similar to the case (1). (3) If $x_1\in D_1$ and $x_2\in D_2$, then $\omega_{x_2}^{D_2}(f^{(2)})=\omega_{x_2}^{D_1\cup D_2}(A)=0$. We can take an $x_3\in D_1\cap D_2\not=\emptyset$ with $\omega_{x_1}^{D_1}\ll\omega_{x_3}^{D_1}$ and $\omega_{x_3}^{D_2}\ll\omega_{x_2}^{D_2}$. Then $\omega_{x_3}^{D_2}(f^{(2)})=0$ which implies $\omega_{x_1}^{D_1\cup D_2}(A)=\omega_{x_1}^{D_1}(f^{(1)})=\omega_{x_3}^{D_1}(f^{(1)})=\omega_{x_3}^{D_2}(f^{(2)})=0$. (4) If $x_1\in D_2$ and $x_2\in D_1$, this case is similar to the case (3).

(b) We will prove $\{\omega_\cdot^{D_1\cup D_2}\}\cup\{\omega_\cdot^{D_\lambda}\}$ is a harmonic measure system.

(b1) We will prove $\omega_x^{D_1\cup D_2}(A)$ is continuous at any $x_0\in D_1\cup D_2$ for any $\omega_\cdot^{D_1\cup D_2}$-measurable $A$. We argue in 2 cases: (1) If $x_0\in D_1$, then there exists a neighborhood $B_{x_0, \delta}\subset D_1$. Hence for any $y\in B_{x_0, \delta}\subset D_1$, we get $\omega_y^{D_1\cup D_2}(A)=\omega_y^{D_1}(f^{(1)})$. Then the continuous of $\omega_y^{D_1}(f^{(1)})$ at $x_0$ implies the continuous of $\omega_y^{D_1\cup D_2}(A)$ at $x_0$. (2) If $x_0\in D_1$, this case is similar to the case (1).

(b2) The proof of the compatibility of $D_1\cup D_2$ is easy from the compatibility of $D_1$ and $D_2$ respectively.
\end{proof}

\vspace{0.2cm}

In the following of this subsection we restrict our discussion in $X=\mathbb{R}^n$. We can construct a new harmonic measure $\omega_\cdot^D$ with $D:=\mathop{\cup}\limits_{k=1}^\infty D_k\subset\mathbb{R}^n$, where $\omega_\cdot^{D_k},k\in\mathbb{N}$ are in a harmonic measure system $\{\omega_\cdot^{D_\lambda}\}\not\ni\omega_\cdot^D$.

\vspace{0.2cm}

Suppose the $\{\omega_\cdot^{D_\lambda}\}$ is a harmonic measure system on $\mathbb{R}^n$ and $\{D_k\}_{k=1}^\infty\subseteq\{D_\lambda\}$ is an increasing sequence with $D:=\cup_{k=1}^\infty D_k\not\in\{D_\lambda\}$ bounded. We extend the harmonic measure $\omega_\cdot^{D_k}$ on $\mathbb{R}^n$, denoted by $\widetilde\omega_\cdot^{D_k}$, such that for any $A\subseteq\mathbb{R}^n$, $\widetilde\omega_\cdot^{D_k}(A):=\omega_\cdot^{D_k}(A\cap\partial D_k)$. Then $\{\widetilde\omega_\cdot^{D_k}\}$ is a sequence of Radon measures on $\mathbb{R}^n$ satisfying 
\begin{equation*}
\sup\limits_k\widetilde\omega_\cdot^{D_k}(K)\leq1<\infty
\end{equation*}
for any compact set $K\subset\mathbb{R}^n$. This implies that for any $x\in D$ there exists a subsequence of $\{\widetilde\omega_x^{D_k}\}$, denoted by $\{\widetilde\omega_x^{D_{k_j^x}}\}$, and a Radon measure $\widetilde\omega^D_x$ on $\mathbb{R}^n$ such that
\begin{equation*}
\widetilde\omega_x^{D_{k_j^x}}\rightharpoonup\widetilde\omega^D_x,
\end{equation*}
that is for any $f\in C_c(\mathbb{R}^n)$,
\begin{equation*}
\lim_{k\rightarrow\infty}\widetilde\omega_x^{D_{k_j^x}}(f)=\widetilde\omega^D_x(f),
\end{equation*}
as $j\rightarrow\infty$. At the same time, by $D_k\nearrow D$ with $D$ bounded and the compatibility of the harmonic measure system, we get for any $f\in C_c(\mathbb{R}^n)$, $\{\widetilde\omega^{D_k}_x(f)\}$ is a Cauchy sequence. Then 
\begin{equation}\label{unifor_c}
\lim_{k\rightarrow\infty}\widetilde\omega^{D_k}_x(f)=\widetilde\omega^{D}_x(f).
\end{equation} 
Define $\omega_\cdot^D:=\widetilde\omega^D_\cdot\,\llcorner\,\partial D$. Since $\omega_\cdot^D$ is a Radon measure, we can get for any $x,y\in D$, $\frac{\mathrm{d}\omega_x^D}{\mathrm{d}\omega_y^D}$ exists, is finite $\omega_y^D$-a.e. and $\omega_y^D$-measurable.

\vspace{0.2cm}

\begin{thm}
Suppose the $\{\omega_\cdot^{D_\lambda}\}$ is a harmonic measure system and $\{D_k\}_{k=1}^\infty\subseteq\{D_\lambda\}$ is an increasing sequence with $D:=\cup_{k=1}^\infty D_k\not\in\{D_\lambda\}$ bounded. Let $\omega_\cdot^D$ be the Radon measure defined as the above. If for any $x,y\in D$, there exists a $C_{x,y}>0$ such that $\frac{\mathrm{d}\omega_x^D}{\mathrm{d}\omega_y^D}(z)\leq C_{x,y}$ for any $z\in\partial D$ and, for any $k_0\in\mathbb{N}$, convergence in \eqref{unifor_c} is uniform with respect to $x\in D_{k_0}$, then $\{\omega_\cdot^D\}\cup\{\omega_\cdot^{D_\lambda}\}$ is a harmonic measure system.
\end{thm}

\vspace{0.2cm}

\begin{proof}
In fact, we need to prove by steps:
\begin{enumerate}
\item[(a)] $\omega^D_\cdot$ is a harmonic measure;
\item[(b)] $\{\omega_\cdot^D\}\cup\{\omega_\cdot^{D_\lambda}\}$ is a harmonic measure system.
\end{enumerate}

\vspace{0.2cm}

(a) $\omega^D_\cdot$ is a harmonic measure. In fact, for any $x\in D$, $\omega_x^D$ is a pre-harmonic measure. 

(a1) Take $f\in C_c(\mathbb{R}^n)$ such that $f(y)=1$ for any $y\in\bar D\backslash D_\varepsilon$ with $D_\varepsilon:=\{y\in D: \mathrm{dist}(y,\partial D)>\varepsilon\}$. Since $D_k\nearrow D$, we get for any $x\in D$ there exists a $k_0\in\mathbb{N}$ such that for any $y_k\in\partial D_k$, $f(y_k)=1$ and $x\in D_k$ for any $k\geq k_0$. Then 
\begin{equation*}
\omega_x^D(\partial D)=\widetilde\omega_x^D(f)=\lim_{k\rightarrow\infty}\widetilde\omega_x^{D_k}(f)=\lim_{k_0\leq k\rightarrow\infty}\widetilde\omega_x^{D_k}(f)=\lim_{k_0\leq k\rightarrow\infty}\omega_x^{D_k}(\partial D_k)=1.
\end{equation*}

(a2) For any $y\in D$ and any $A\subset\partial D$, if $\omega^D_y(A)>0$, then taking $0\leq f\in C_c(\mathbb{R}^n)$ with $A\subseteq\mathrm{supp} f$, we get
\begin{equation*}
0<\varepsilon_{f,D}=:\omega^D_y(f)=\widetilde\omega^D_y(f)=\lim_{k\rightarrow\infty}\widetilde\omega_y^{D_k}(f).
\end{equation*}
Since $D_k\nearrow D$ and $\omega^{D_k}_y(f)\rightarrow\omega_y^D(f)$, for any $x,y\in D$ there exists a $k_0\in\mathbb{N}$, such that for any $k\geq k_0$, $x,y\in D_k$ and
\begin{equation*}
\omega_y^{D_k}(f)=\widetilde\omega_y^{D_k}(f)\geq\frac{\varepsilon_{f,D}}{2}>0.
\end{equation*}
By the definition of $\omega_\cdot^D$, $\frac{\mathrm{d}\omega_y^{D_k}}{\mathrm{d}\omega_x^{D_k}}(z_k)\rightarrow\frac{\mathrm{d}\omega_x^D}{\mathrm{d}\omega_y^D}(z)$ as $k\rightarrow\infty$. Since $\frac{\mathrm{d}\omega_x^D}{\mathrm{d}\omega_y^D}$ is bounded on $\partial D$, there exists a constant $\widetilde C_{x,y}>0$ such that for any $k\geq k_0$ and any $z_k\in\partial D_k$,
\begin{equation*}
\frac{\mathrm{d}\omega_y^{D_k}}{\mathrm{d}\omega_x^{D_k}}(z_k)\leq\widetilde C_{x,y}. 
\end{equation*}
At the same time since $\omega_x^{D_k}\ll\omega_y^{D_k}$, we get for any $k\geq k _0$,
\begin{equation*}
\omega_x^{D_k}(f)=\int_{\partial D_k}f\,\mathrm{d}\omega_x^{D_k}=\int_{\partial D_k}f\frac{\mathrm{d}\omega_x^{D_k}}{\mathrm{d}\omega_y^{D_k}}\,\mathrm{d}\omega_y^{D_k}\geq \frac{1}{\widetilde C_{x,y}}\omega_y^{D_k}(f)\geq \frac{1}{\widetilde C_{x,y}}\frac{\varepsilon_{f,D}}{2}.
\end{equation*}
Sending $k\rightarrow\infty$ we get 
\begin{equation*}
\omega_x^{D}(f)=\lim_{k\rightarrow\infty}\omega_x^{D_k}(f)\geq\frac{1}{\widetilde C_{x,y}}\frac{\varepsilon_{f,D}}{2}>0,
\end{equation*}
which implies $\omega_x^D(A)>0$. Then we get $\omega_y^{D}\ll\omega_x^{D}$ for any $x,y\in D$.

\vspace{0.2cm}

(b) $\{\omega_\cdot^D\}\cup\{\omega_\cdot^{D_\lambda}\}$ is a harmonic measure system. 

(b1) We will prove for any $x_0\in D$ and $f\in C_c(\mathbb{R}^n)$, $\omega_{x}(f)$ is continuous at $x_0$. In fact, by $D_k\nearrow D$ and the uniform convergence of \eqref{unifor_c}, for any $\varepsilon>0$, there exits a $k_0\in\mathbb{N}$, such that $B_{x_0, \frac{1}{2}\mathrm{dist}(x_0,\partial D_{k_0})}\subseteq D_{k_0}$,
\begin{equation*}
|\omega_{x_0}^D(f)-\omega_{x_0}^{D_{k_0}}(f)|<\frac{\varepsilon}{3},
\end{equation*}
and for any $y\in B_{x_0, \frac{1}{2}\mathrm{dist}(x_0,\partial D_{k_0})}$
\begin{equation*}
|\omega_{y}^D(f)-\omega_{y}^{D_{k_0}}(f)|<\frac{\varepsilon}{3}.
\end{equation*}
For the above $\varepsilon$ and $k_0$, since $\omega_x^{D_{k_0}}$ is continuous at $x_0$, we get there exists a $\delta=\delta(k_0,x_0,\varepsilon)>0$ such that for any $y\in B_{x_0,\delta}$
\begin{equation*}
|\omega_{x_0}^{D_{k_0}}(f)-\omega_{y}^{D_{k_0}}(f)|<\frac{\varepsilon}{3}.
\end{equation*}
Then for the above $\varepsilon>0$, take $\delta_0:=\min\{ \frac{1}{2}\mathrm{dist}(x_0,\partial D_{k_0}),\delta\}$, for any $y\in B_{x_0,\delta_0}$ we get
\begin{equation*}
\begin{aligned}
&|\omega_{x_0}^D(f)-\omega_y^D(f)|\\
\leq&|\omega_{x_0}^D(f)-\omega_{x_0}^{D_{k_0}}(f)|+|\omega_{x_0}^{D_{k_0}}(f)-\omega_{y}^{D_{k_0}}(f)|+|\omega_{y}^{D_{k_0}}(f)-\omega_{y}^{D}(f)|\\
\leq&\frac{\varepsilon}{3}+\frac{\varepsilon}{3}+\frac{\varepsilon}{3}=\varepsilon.
\end{aligned}
\end{equation*}

(b2) The proof of the compatibility of $D$ is easy from the compatibility of $D_k$.
\end{proof}

\vspace{0.2cm}

\subsection{Harnack inequality and Harnack principle}\label{sec2.4}

\vspace{0.2cm}

In this subsection we consider Harnack inequality and Harnack principle in a harmonic measure system.

\vspace{0.2cm}

\begin{defi}
Given a harmonic measure system $\{\omega_\cdot^{D_\lambda}\}$ on the Hausdorff space $X$. Suppose $D\in \{D_\lambda\}$ with an subset $A\subseteq D$. $H:A\subseteq D\mapsto H(A;D)\in[1,+\infty]$ defined by
\begin{equation*}
H(A;D):=\sup_{x,y\in A\atop 0\leq f\in L(\partial D;\omega_\cdot^D)\backslash\{0\}}\frac{\omega_x^D(f)}{\omega_y^D(f)}
\end{equation*}
is called the \emph{Harnack index} of the subset $A$ relative to $D$. Denote
\begin{equation*}
D_\alpha:=\{A\subseteq D:\,H(A;D)=\alpha\}.
\end{equation*}
The elements of $D_\alpha$ are called the \emph{Harnack subsets} with index $\alpha$ relative to $D$.
\end{defi}

\vspace{0.2cm}

By the definition of Harnack index, the following propositions holds.

\vspace{0.2cm}

\begin{prop}\label{Harnack_basic}
\begin{enumerate}
\item[(1)] For any $D\in\{D_\lambda\}$, $H(D;D)=+\infty$;
\item[(2)] If $D_1,D_2\in\{D_\lambda\}$ with a subset $A\subseteq D_2\subseteq D_1$, then $H(A;D_1)\leq H(A;D_2)$ and $H(A;D_1)\leq H(D_2;D_1)$;
\item[(3)] For any $D\in\{D_\lambda\}$, any $f\in L(\partial D,\omega_\cdot^D)\backslash\{0\}$ with $f\geq0$ and any $x,y\in A\subseteq D$,
\begin{equation*}
\omega_x^{D}(f)\leq H(A;D)\omega_y^{D}(f).
\end{equation*}
\end{enumerate}
\end{prop}

\vspace{0.2cm}

Suppose $\{\omega^{D_\lambda}_\cdot\}$ is a Radon harmonic measure system on $X$. Radon-Nikodym Theorem implies $\frac{\mathrm{d}\omega^D_x}{\mathrm{d}\omega^D_y}$ exists. The boundedness of $\frac{\mathrm{d}\omega^D_x}{\mathrm{d}\omega^D_y}$ implies $H(A;D)<+\infty$.

\vspace{0.2cm}

\begin{thm}[General Harnack inequality]\label{Ha}
Suppose $D\in\{D_\lambda\}$. If $\frac{\mathrm{d}\omega^D_x}{\mathrm{d}\omega^D_y}(z)$ is bounded with respect to $x,y\in A\subseteq D$ and $z\in\partial D$, then $H(A;D)<+\infty$.
\end{thm}

\vspace{0.2cm}

\begin{proof}
By the boundedness of $\frac{\mathrm{d}\omega^D_x}{\mathrm{d}\omega^D_y}$, for any $f\in L(\partial D,\omega_\cdot^D)\backslash\{0\}$ with $f\geq0$ and any $x,y\in A\subseteq D$, Radon-Nikodym Theorem implies
\begin{equation*}
\frac{\omega^D_x(f)}{\omega^D_y(f)}=\frac{\int f\,\mathrm{d}\omega_x^D}{\int f\,\mathrm{d}\omega_y^D}=\frac{\int f\frac{\mathrm{d}\omega^D_x}{\mathrm{d}\omega^D_y}\,\mathrm{d}\omega_y^D}{\int f\,\mathrm{d}\omega_y^D}\leq\frac{\int Mf\,\mathrm{d}\omega_y^D}{\int f\,\mathrm{d}\omega_y^D}=M.
\end{equation*}
Then $H(A;D)\leq M<\infty$.
\end{proof}

\vspace{0.2cm}

\begin{cor}[Harnack inequality]
Suppose $D\in\{D_\lambda\}$ and $\partial D$ is compact. If $\frac{\mathrm{d}\omega^D_x}{\mathrm{d}\omega^D_y}(z)$ is continuous with respect to $x,y\in D$ and $z\in\partial D$. Then $H(A;D)<+\infty$ for any subset $A\subset\subset D$.
\end{cor}
\vspace{0.2cm}

In the following of this subsection, we restrict our discussion in $X=\mathbb{R}^n$.

\vspace{0.2cm}

\begin{enumerate}
\item[(TSCI)] (Translation and Scaling Closed \& Invariant) The harmonic measure system $\{\omega_\cdot^{D_\lambda}\}$ in $\mathbb{R}^n$ is closed and invariant under the transformation of translation and scaling, that is, for any $x_0,\nu_0\in\mathbb{R}^n$ and $k\in\mathbb{R} (k\not=0)$, define the translation and scaling map for any $x\in\mathbb{R}^n$ by $\nu_{\nu_0,x_0,k} (x):=\nu_0+k(x-x_0)$, then 
\begin{enumerate}
\item[(1)] (Translation and Scaling Closed) For any open $D\in\{D_\lambda\}$, $\nu_{\nu_0,x_0,k}(D)\in\{D_\lambda\}$; 
\item[(2)] (Translation and Scaling Invariant) For any $\omega_\cdot^D$-measurable set $A\subseteq\partial D$, $\nu_{\nu_0,x_0,k}(A)$ is $\omega_\cdot^{\nu_{\nu_0,x_0,k}(D)}$-measurable with
\begin{equation*}
\omega_{\nu_{\nu_0,x_0,k}(x)}^{\nu_{\nu_0,x_0,k}(D)}(\nu_{\nu_0,x_0,k}(A)) = \omega_{x}^{D}(A).
\end{equation*}
\end{enumerate}
\end{enumerate}

\vspace{0.2cm}

Roughly speaking, the TSCI assumption on the harmonic measure system means that the harmonic measure $\omega_{\cdot}^D$ in the harmonic measure system depends on the shape of $D\in\{D_\lambda\}$ but independent on its size and location in $\mathbb{R}^n$.

\vspace{0.2cm}

\begin{thm}[Harnack principle in $\mathbb{R}^n$]
Suppose the Radon harmonic measure system $\{\omega_\cdot^{D_\lambda}\}$ in $\mathbb{R}^n$ satisfies TSCI and $\frac{\mathrm{d}\omega^D_x}{\mathrm{d}\omega^D_y}(z)$ is continuous with respect to $x,y\in D\in\{D_\lambda\}$ and $z\in\partial D$. If $D,D_1\in\{D_\lambda\}$ with a subset $A\subset\subset D_1\subseteq D$ and $\partial D$ compact, then $H(\nu(A);D)\leq H(A;D_1)<+\infty$ provided $\nu(D_1)\subseteq D$ for any translation and scaling map $\nu$. Equivalently, for any $f\in L(\partial D;\omega_\cdot^D)$ and any $x,y\in\nu(A)$, there exists a constant $C:=H(A;D_1)$ such that
\begin{equation*}
\omega_x^D(f)\leq C\omega_y^D(f).
\end{equation*}
\end{thm}

\vspace{0.2cm}

\begin{proof}
The TSCI assumption on the harmonic measure system $\{\omega_\cdot^{D}\}$, Proposition \ref{Harnack_basic} (2) and Theorem \ref{Ha} imply
\begin{equation*}
H(\nu(A);D)\leq H(\nu(A);\nu(D_1))=H(A;D_1)<+\infty.
\end{equation*}
\end{proof}

\vspace{0.2cm}

\begin{rem}
In the applications when we use the Harnack principle in $\mathbb{R}^n$, we always say Harnack inequality holds for the constant $C$ depending on the shape of $A\subset\subset D_1$. See a special inner cone case of the Harnack principle in Proposition \ref{Harnack_Lip}.
\end{rem}

\vspace{0.2cm}

\section{Harmonic measure system induced by the elliptic operator}\label{Sec3}

\vspace{0.2cm}

In this section, we give some examples of the harmonic measure system induced by the elliptic operator.

\vspace{0.2cm}

\subsection{The classical Laplace operator in $\mathbb{R}^n$}

\vspace{0.2cm}

The most simple and important one is the harmonic measure system induced by classical Laplace operator $-\Delta:=-\mathrm{div}\nabla$ in $\mathbb{R}^n$. We continue to use the notations in Subsection \ref{sec1_1}.

\vspace{0.2cm}

Suppose $D$ is a bounded open set in $\mathbb{R}^n$. According to Theorem \ref{Wiener}, each $f\in C(\partial D)$ determines a $H_f(x)$, where $x\in D$. Consider the mapping $L_x:C(\partial D)\rightarrow\mathbb{R}$ defined by $L_x:=H_f(x)$. By the Reisz representation theorem for measures, we can get a Borel measure $\omega_x^D$ on $\partial D$ in the following lemma.

\vspace{0.2cm}

\begin{lem}[Lemma 8.12 in \cite{Helm69}]\label{lem31}
For any $x\in D$, $L_x$ is a positive linear functional on $C(\partial D)$; Moreover, there is a unique unit measure $\omega_x^D$ on the Borel subsets of $\partial D$ such that $H_f(x)=L_x(f)=\omega_x^D(f)$ for any $f\in C(\partial D)$.
\end{lem}

\vspace{0.2cm}

We list some propositions related to the Borel measure $\omega_x^D$ on $\partial D$ in Lemma \ref{lem31}.

\vspace{0.2cm}

\begin{thm}[Theorem 8.13 in \cite{Helm69}]
If $f$ is a lower bounded (or upper bounded) Borel measurable function on $\partial D$, then $\overline H_f(x)=\underline H_f(x)=\omega_x^D(f)$ for any $x\in D$.
\end{thm}

\vspace{0.2cm}

\begin{thm}[Theorem 8.14 in \cite{Helm69}]\label{lem33}
If $V$ is a component of the bounded open set $D$, then the class of Borel subsets of $\partial D$ of $\omega_x^D$-measure zero is independent of $x\in V$.
\end{thm}

\vspace{0.2cm}

\begin{thm}[Brelot, Theorem 8.17 in \cite{Helm69}]\label{th33}
A boundary function $f$ is resolutive if and only if it is $\omega_\cdot^D$-integrable; in which case $H_f(x)=\omega_x^D(f)$ for any $x\in D$.
\end{thm}

\vspace{0.2cm}

\begin{lem}[Corollary 8.21 in \cite{Helm69}]\label{lem35}
Suppose there is a barrier at $x\in\partial D$. If $f$ is bounded on $\partial D$ and continuous at $x$, then
\begin{equation*}
\lim_{D\ni y\rightarrow x}\overline H_f(y)=\lim_{D\ni y\rightarrow x}\underline H_f(y)=f(x).
\end{equation*}
\end{lem}

\vspace{0.2cm}

\begin{thm}[Theorem 8.22 in \cite{Helm69}]\label{lem36}
A point $x\in\partial D$ is a regular boundary point in the sense of Definition \ref{H_m} if and only if there is a barrier at $x$.
\end{thm}

\vspace{0.2cm}

\begin{lem}[Corollary 8.28 in \cite{Helm69}]\label{lem37}
If $D$ is an open subset of $\mathbb{R}^n$, $n\geq2$, then there is an increasing sequence $\{D_j\}$ of bounded regular (in the sense of Definition \ref{H_m}) open sets with closures in $D$ such that $D=\cup D_j$.
\end{lem}

\vspace{0.2cm}

By Lemma \ref{lem31} and Theorem \ref{lem33}, we get $\omega_x^D$ is a pre-harmonic measure in the sense of Definition \ref{pre_HM}. And $\omega_\cdot^D$ in Theorem \ref{th33} is a harmonic measure in the sense of Definition \ref{H_m}. The harmonic measure system in this case is 
\begin{equation*}
\{\omega_\cdot^D: D \text{ is any bounded open set of }\mathbb{R}^n\}.
\end{equation*} 
Theorem \ref{lem36} insures bounded Lipschitz domain is regular in the sense of Definition \ref{H_m}. By Lemma \ref{lem37}, the harmonic measure system $\{\omega_\cdot^D\}$ is regular closed.

\subsection{A class of degenerate elliptic operator in $\mathbb{R}^n$}

\vspace{0.2cm}

In this subsection, a harmonic measure on sets with co-dimension higher than $1$ is considered, see \cite{Davi21}. The harmonic measure is induced by a kind of degenerate elliptic operator.

\vspace{0.2cm}

A set $\Gamma\subset\mathbb{R}^n$ is a $d$-dimensional \emph{Ahlfors Regular} (AR) if there exists a measure $\sigma$ supported on $\Gamma$ and a constant $C_0\geq 1$ such that
\begin{equation*}
C_0^{-1}r^d\leq\sigma(B(x,r))\leq C_0r^d
\end{equation*}
for $x\in\Gamma$ and $0<r<\mathrm{diam}(\Gamma)$.

\vspace{0.2cm}

Suppose $\Gamma$ is an AR set of dimension $d<n-1$ and $D:=\mathbb{R}^n\backslash\Gamma$. Define a weight $w(x):=\dist(x,\partial D)^{d+1-n}$ and a measure $m(E):=\int_Ew(x)\,\mathrm{d}x$. Define the weighted Sobolev space $W$ by
\begin{equation*}
W:=\{u\in L_{\mathrm{loc}}^1:\nabla u\in L^2(D,\mathrm{d}m)\},
\end{equation*}
and the set of traces $H$, which is the set of measurable functions $g$ defined on $\Gamma$ such that
\begin{equation*}
\|g\|_H^2:=\int_\Gamma\int_\Gamma\frac{|g(x)-g(y)|^2}{|x-y|^{d+1}}\,\mathrm{d}\sigma(x)\,\mathrm{d}\sigma(y)<+\infty.
\end{equation*}
For $\sigma$-almost every $x\in\Gamma$, define the trace operator $\mathrm{Tr}:W\rightarrow H$ at $x$ by 
\begin{equation*}
\mathrm{Tr}\,u(x):=\lim_{r\rightarrow0}\fint_{B(x,r)}u(y)\,\mathrm{d}y
\end{equation*}

\vspace{0.2cm}

Now consider any degenerate elliptic operator $L:=-\mathrm{div} A\nabla$ with measurable real matrix-valued coefficients $A$ satisfying there exists a constant $C_1\geq 1$ such that
\begin{equation*}
A(x)\xi\cdot\eta\leq C_1w(x)|\xi||\eta|, \,\forall\, x\in D, \forall\,\xi,\eta\in\mathbb{R}^n
\end{equation*}
and
\begin{equation*}
A(x)\xi\cdot\xi\geq C_1^{-1}w(x)|\xi|^2, \,\forall\, x\in D, \forall\,\xi\in\mathbb{R}^n.
\end{equation*}
By the Lax-Milgram theorem, for each $g\in H$, one can find a unique $u\in W$ such that
\begin{equation}\label{weak_s}
\left\{\begin{aligned}
&Lu=0&&\text{ in } D,\\
&\mathrm{Tr}\,u=g&&\text{ on } \Gamma,
\end{aligned}\right.
\end{equation}
where $Lu=0$ is in the sense of \emph{weak solution}:
\begin{equation*}
\int_DA\nabla u\cdot\nabla\varphi\,\mathrm{d}x=0, \,\forall\,\varphi\in C_0^\infty(D).
\end{equation*}

\vspace{0.2cm}

In \cite{Davi21} the authors get the following regularity result for the weak solution for \eqref{weak_s}.

\vspace{0.2cm}

\begin{thm}[De Giorgi-Nash-Moser estimates at the boundary, see \cite{Davi21}]\label{t35}
Let $B=B(x,r)$ be a ball centered on $\Gamma$ and $u\in W\subset L_{\mathrm{loc}}^2(\mathbb{R}^n,\mathrm{d}m)$ be a weak solution to $Lu=0$ in $D$ such that $\mathrm{Tr}\,u$ is continuous and bounded on $B$. Denote by $\mathrm{osc}_Eu$ the difference between the (essential) supremum and the (essential) infimum of $u$ on $E$. There exists $\alpha>0$ such that for $0<s<r$,
\begin{equation*}
\mathop{\mathrm{osc}}\limits_{B(x,s)}u\leq C\left(\frac{s}{r}\right)^\alpha\mathop{\mathrm{osc}}\limits_{B(x,r)}u+C\mathop{\mathrm{osc}}\limits_{B(x,\sqrt{sr})\cap\Gamma}\mathrm{Tr}\,u.
\end{equation*}
In particular, $u$ is continuous on $B$. In addition, if $\mathrm{Tr}\,u=0$ on $B$, then for $x\in\frac{1}{2}B$ and $0<s<\frac{r}{3}$,
\begin{equation*}
\mathop{\mathrm{osc}}\limits_{B(x,s)}u\leq C\left(\frac{s}{r}\right)^\alpha\left(\fint_B|u|^2\mathrm{d}m\right)^{\frac{1}{2}}<+\infty.
\end{equation*}
The constants $\alpha, C$ depend on the dimensions $d$ and $n$, as well as the constants $C_0$ and $C_1$.
\end{thm}

\vspace{0.2cm}

For any $x\in D$, the unique solvability of \eqref{weak_s} and Theorem \ref{t35} determine a mapping $U:H\cap C_c(\Gamma)\rightarrow C(\mathbb{R}^n)$, where $C_c(\Gamma)$ denotes the space of continuous functions with compact support on $\Gamma$, and $Ug$ is the unique solution to \eqref{weak_s}. The following lemma extend the mapping $U$ to $C_c(\Gamma)$.

\vspace{0.2cm}

\begin{lem}[Lemma 9.23 in \cite{Davi21}]\label{lem_3_4}
There exists a bounded linear operator
\begin{equation*}
U: C_c(\Gamma)\rightarrow C(\mathbb{R}^n)
\end{equation*}
such that, for any $g\in C_c(\Gamma)$,
\begin{enumerate}
\item[(1)] $(Ug)|_\Gamma=g$;
\item[(2)] $\sup\limits_{\mathbb{R}^n}Ug=\sup\limits_\Gamma g$ and $\inf\limits_{\mathbb{R}^n}Ug=\inf\limits_\Gamma g$;
\item[(3)] $Ug\in W_{\mathrm{loc}}(D)$ and is a weak solution of $Lu=0$ in $D$;
\item[(4)] If $B$ is a ball centered on $\Gamma$ and $g=0$ on $B$, then $Ug\in W_{\mathrm{loc}}(B)$;
\item[(5)] If $g\in H\cap C_c(\Gamma)$, then $Ug\in W$.
\end{enumerate}
\end{lem}

\vspace{0.2cm}

Now for any $x\in D$, consider the mapping $L_x:C_c(\Gamma)\rightarrow\mathbb{R}$ defined by $L_x(g):=Ug(x)$. By the Reisz representation theorem for measures, we can get a Borel measure $\omega_x^D$ on $\Gamma$ such that
\begin{equation*}
Ug(x)=\omega_x^D(g)
\end{equation*}
for any $g\in C_c(\Gamma)$ (Lemma 9.30 in \cite{Davi21}).

\vspace{0.2cm}

We list some properties of the Borel measure $\omega_x^D$.

\vspace{0.2cm}

\begin{lem}[Lemma 9.33 and Lemma 9.38 in \cite{Davi21}]\label{lem3_5}
(1) For any $x\in D$, $\omega^D_x(\Gamma)=1$; (2) For any Borel set $E\subset\Gamma$, if there exists an $x\in D$ such that $\omega_x^D(E)=0$, then for any $y\in D$, $\omega_y^D(E)=0$.
\end{lem}

\vspace{0.2cm}

\begin{lem}[Lemma 9.38 in \cite{Davi21}]
Let $E\subset\Gamma$ be a Borel set. Then
\begin{enumerate}
\item[(1)] $\omega^D_x(E)\in W_{\mathrm{loc}}(D)$ and is a weak solution of $Lu=0$ in $D$; 
\item[(2)] If $B\subset\mathbb{R}^n$ is a ball such that $E\cap B=\emptyset$, then $\omega^D_x(E)\in W_{\mathrm{loc}}(B)$ and $\mathrm{Tr}\,\omega_{x}^D(E)=0$ on $B\cap\Gamma$. 
\end{enumerate}
\end{lem}

\vspace{0.2cm}

Lemma \ref{lem3_5} indicates $\omega^D_x$ is a pre-harmonic measure in the sense of Definition \ref{pre_HM}. By the theory developed in Section \ref{Sec2}, we can construct a harmonic measure system. And the corresponding maximum principle, Harnack inequality and Harnack principle hold in the induced harmonic measure system.

\vspace{0.2cm}

\subsection{Dirichlet boundary problem in graph theory}

\vspace{0.2cm}

Let $X$ be a finite set, whose elements are called \emph{vertexes}. Define a \emph{weight function} $\mu:X\times X\rightarrow[0,\infty), (x,y)\mapsto\mu_{xy}$ which is symmetric, has zero diagonal. $\mu_{xy}>0$ represents there is an \emph{edge} between $x$ and $y$ and $\mu_{xy}=0$ represents there is no edge between $x$ and $y$. $X$ with the edges induced by $\mu$ is a \emph{edge weighted graph}, denoted by $(X,\mu)$. If $\mu_{xy}>0$, we write $x\sim y$ and let $xy$ and $yx$ be the \emph{oriented edges} of the graph. We assume $(X,\mu)$ is connected, that is for any two vertexes $x,y\in X$, there is a path $x=x_0\sim x_1\sim x_2\sim\cdots\sim x_n=y$. 

\vspace{0.2cm}

It is well known that $(X,\mu)$ can be embedded on a compact surface in an Euclidean space. We define the topology of $(X,\mu)$ be the induced topology by the embedded surface. For any connected subgraph $(X',\mu':=\mu|_{X'\times X'})$ of $(X,\mu)$, denote $\bar X':=X'\cup B$, where $B:=\{y\in X\backslash X':y\sim x\in X'\}$ is called the \emph{relative boundary} of $(X',\mu')$. $(\bar X',\bar\mu':=\mu|_{\bar X'\times\bar X'})^\circ:=(\bar X',\bar\mu')\backslash(B,\mu_B:=\mu|_{B\times B})$ is called \emph{open} subgraph of $(X,\mu)$ containing $(X',\mu')$. Under the topology induced by the topology of $(X,\mu)$, $\partial(\bar X',\bar\mu')^\circ=B$.

\vspace{0.2cm}

Suppose $(X',\mu')$ is a connected subgraph of $(X,\mu)$. We define the \emph{Laplacian} $\Delta_{(\bar X',\bar\mu')}$ on the set $U(\bar X'):=\{u:\bar X'\rightarrow\mathbb{R}\}$ for any $x\in(\bar X',\bar\mu')^{\circ}\cap\bar X'=X'$ by
\begin{equation*}
\Delta_{(\bar X',\bar\mu')} u(x)=\sum_{x\sim y\in\bar X'}\bar\mu'(x,y)(u(x)-u(y)).
\end{equation*}
A function $u:\bar X'\rightarrow\mathbb{R}$ is called \emph{harmonic} at $x\in X'$ if $u\in U(\bar X')$ satisfies $\Delta_{(\bar X',\bar\mu')}u(x)=0$.

\vspace{0.2cm}

The \emph{Dirichlet boundary problem} on $(\bar X',\bar\mu')$ is, for a given $f:B\rightarrow\mathbb{R}$, to find a function $u\in U(\bar X')$ such that
\begin{equation*}
(P)\left\{\begin{aligned}
&\Delta_{(\bar X',\bar\mu')} u(x)=0 &&\text{ for } x\in X';\\
&u(x)=f(x) &&\text{ for } x\in B.
\end{aligned}\right.
\end{equation*}

\vspace{0.2cm}

\begin{lem}\label{Harp}
For any $x\in (\bar X',\bar\mu')^{\circ}\cap \bar X'$, the solution $u$ of $(P)$ satisfies
\begin{equation*}
u(x)=\sum_{x\sim y\in\bar X'}\frac{\bar\mu'(x,y)}{\sum\limits_{x\sim z\in\bar X'}\bar\mu'(x,z)}u(y)
\end{equation*}
\end{lem}

\vspace{0.2cm}

\begin{thm}\label{EU}
If $(\bar X',\bar\mu')$ is connected and $\partial(\bar X',\bar\mu')\not=\emptyset$, then the solution to $(P)$ exists uniquely.
\end{thm}

\vspace{0.2cm}

\begin{proof}
Denote $(\bar X',\bar\mu')^{\circ}\cap \bar X'=\{x_i\}_{i=1}^n$, $\bar\mu'(x_i,x_j):=a_{ij}$. By Lemma \ref{Harp}, the existence and uniqueness of the solution to $(P)$ is equivalent to the existence and uniqueness of the solution to linear system $Ax=y$, where 
\begin{equation*}
A:=\begin{bmatrix}
k_1+\sum_{j=1}^n a_{1j} &-a_{12}&-a_{13}&\cdots &-a_{1n}\\
-a_{21}&k_2+\sum_{j=1}^n a_{2j}&-a_{23}&\cdots &-a_{2n}\\
-a_{31}&-a_{32}&k_3+\sum_{j=1}^n a_{3j}&\cdots &-a_{3n}\\
\vdots&\vdots&\vdots&\ddots&\vdots\\
 -a_{n1}&-a_{n2}&-a_{n3}&\cdots &k_n+\sum_{j=1}^n a_{nj}
\end{bmatrix},
\end{equation*}
\begin{equation*}
k_i:=\sum_{x_j\sim y\in\partial(\bar X',\bar\mu')}\bar\mu'(x_i,y),
\end{equation*}
\begin{equation*}
x:=\begin{bmatrix}
u(x_1)\\u(x_2)\\\vdots\\u(x_n)
\end{bmatrix},
\end{equation*}
\begin{equation*}
y:=\begin{bmatrix}
\sum\limits_{x_1\sim y\in\partial(\bar X',\bar\mu')}\bar\mu'(x_1,y)f(y)\\
\sum\limits_{x_2\sim y\in\partial(\bar X',\bar\mu')}\bar\mu'(x_2,y)f(y)\\
\vdots\\
\sum\limits_{x_n\sim y\in\partial(\bar X',\bar\mu')}\bar\mu'(x_n,y)f(y)
\end{bmatrix}.
\end{equation*}
It is clear that $k_i\geq0$, $a_{ij}\geq0$.

We will prove $A$ is positive, that is, for any $0\not=\xi=(\xi_1,\xi_2,\cdots,\xi_n)\in\mathbb{R}^n$, $\xi A\xi^T>0$. In fact,
\begin{equation*}
\xi A\xi^T=\sum_{i=1}^nk_i\xi_i^2+\sum_{i,j=1}^na_{ij}(\xi_i-\xi_j)^2.
\end{equation*}

(1) If $\sum\limits_{i,j=1}^na_{ij}(\xi_i-\xi_j)^2>0$, then $A$ is positive. 

(2) if $\sum\limits_{i,j=1}^na_{ij}(\xi_i-\xi_j)^2=0$, then the connectivity of $(\bar X',\bar\mu')$ implies that $\xi_1=\xi_2=\cdots=\xi_n\not=0$. By $\partial(\bar X',\bar\mu')\not=\emptyset$, we know there exists a $k_l\not=0$. Then $\sum\limits_{i=1}^nk_i\xi_i^2>0$, which implies that $A$ is positive.

By (1) and (2), $A$ is positive. Then $|A|\not=0$, which implies the existence and uniqueness of the solution to $Ax=y$, equivalently, solution to $(P)$ exists uniquely.
\end{proof}

\vspace{0.2cm}

For a simple graph $(X,\mu)$ with $X=\{x_1,x_2\}$ and $\mu(x_1,x_2)>0$, we add a new vertex $x_\lambda\in x_1x_2$ in $X$ in the following way:
\begin{equation}\label{Cd1}
\frac{\widetilde\mu(x_1,x_\lambda)}{\widetilde\mu(x_\lambda,x_2)}=\frac{1-\lambda}{\lambda}:=\frac{|x_\lambda x_2|}{|x_1x_\lambda|}.
\end{equation}
\begin{equation}\label{Cd2}
\frac{1}{\widetilde\mu(x_1,x_\lambda)}+\frac{1}{\widetilde\mu(x_\lambda,x_2)}=\frac{1}{\mu(x_1,x_2)},
\end{equation}
Then $(\widetilde X:=\{x_1,x_\lambda,x_2\},\widetilde\mu)$ is a new edge weighted graph. \eqref{Cd1} ensures the following lemma holds.

\vspace{0.2cm}

\begin{lem}\label{Addp}
In the graph $(\widetilde X,\widetilde\mu)$, $u(x_\lambda):=(1-\lambda)u(x_1)+\lambda u(x_2)$ is harmonic at $x_\lambda$, where $\lambda$ is defined in \eqref{Cd1}.
\end{lem}

\vspace{0.2cm}

\begin{proof}
By the Lemma \ref{Harp} and \eqref{Cd1}, if $\Delta_{(\widetilde X,\widetilde\mu)}u(x_\lambda)=0$, then
\begin{equation*}
\begin{aligned}
u(x_\lambda)&=\frac{\widetilde\mu(x_\lambda,x_1)}{\widetilde\mu(x_\lambda,x_1)+\widetilde\mu(x_\lambda,x_2)}u(x_1)+\frac{\widetilde\mu(x_\lambda,x_2)}{\widetilde\mu(x_\lambda,x_1)+\widetilde\mu(x_\lambda,x_2)}u(x_2)\\
&=(1-\lambda)u(x_1)+\lambda u(x_2).
\end{aligned}
\end{equation*}
\end{proof}

\vspace{0.2cm}

\eqref{Cd1} and \eqref{Cd2} ensures the following lemma holds.

\vspace{0.2cm}

\begin{lem}\label{bdharmo}
Suppose $x_1x_2$ is an edge in the subgraph $(X',\mu')$ of $(X,\mu)$, $x_1\in(\bar X',\bar\mu')^\circ$ and $u$ is harmonic at $x_1$. If $x_\lambda\in x_1x_2$ and $u(x_\lambda):=(1-\lambda)u(x_1)+\lambda u(x_2)$, then $u$ is harmonic at $x_1$ in $(\overline{X'\cup\{x_\lambda\}},\widetilde{\bar\mu'})$, where $\widetilde{\bar\mu'}(x_1,x_\lambda)$ and $\widetilde{\bar\mu'}(x_2,x_\lambda)$ satisfies conditions as \eqref{Cd1} and \eqref{Cd2}, on the other vertexes $\widetilde{\bar\mu'}=\bar\mu'$.
\end{lem}

\vspace{0.2cm}

\begin{proof}
Since $u$ is harmonic at $x_1$ in $(\bar X',\bar\mu')$, 
\begin{equation*}
\begin{aligned}
0&=\Delta_{(\bar X',\bar\mu')} u(x_1)=\sum_{x_1\sim y\in\bar X'}\bar\mu'(x_1,y)(u(x_1)-u(y))\\
&=\bar\mu'(x_1,x_2)(u(x_1)-u(x_2))+\sum_{x_1\sim y\not=x_2,y\in\bar X'}\bar\mu'(x_1,y)(u(x_1)-u(y)).
\end{aligned}
\end{equation*}
In $(\overline{X'\cup\{x_\lambda\}},\widetilde{\bar\mu'})$, we get
\begin{equation*}
\begin{aligned}
\Delta_{(\overline{X'\cup\{x_\lambda\}},\widetilde{\bar\mu'})} u(x_1)&=\sum_{x_1\sim y\in\overline{X'\cup\{x_\lambda\}}}\widetilde{\bar\mu'}(x_1,y)(u(x_1)-u(y))\\
&=\widetilde{\bar\mu'}(x_1,x_\lambda)(u(x_1)-u(x_\lambda))\\
&\quad\quad+\sum_{x_1\sim y\not=x_\lambda,y\in\overline{X'\cup\{x_\lambda\}}}\widetilde{\bar\mu'}(x_1,y)(u(x_1)-u(y))\\
&=A(x_1,x_\lambda,x_2)\bar\mu'(x_1,x_2)(u(x_1)-u(x_2))\\
&\quad\quad+\sum_{x_1\sim y\not=x_2,y\in\bar X'}\bar\mu'(x_1,y)(u(x_1)-u(y)),
\end{aligned}
\end{equation*}
where 
\begin{equation*}
A(x_1,x_\lambda,x_2):=\frac{\widetilde{\bar\mu'}(x_1,x_\lambda)}{\bar\mu'(x_1,x_2)}\frac{u(x_1)-u(x_\lambda)}{u(x_1)-u(x_2)}.
\end{equation*}
By \eqref{Cd1} and \eqref{Cd2}, it is easy to get $A(x_1,x_\lambda,x_2)=\dfrac{1}{\lambda}\cdot\lambda=1$. Then $\Delta_{(\overline{X'\cup\{x_\lambda\}},\widetilde{\bar\mu'})} u(x_1)=\Delta_{(\bar X',\bar\mu')} u(x_1)=0$, which completes the proof.
\end{proof}

In general, we can extend the harmonic function on vertexes $X'$ to points on edges of $(\bar X',\bar\mu')^{\circ}$ by Lemma \ref{Addp} and Lemma \ref{bdharmo}. Then for any open connected set $O$ of the compact surface the graph $(X,\mu)$ embedded on, and any $x\in D:=(X,\mu)\cap O$ with $D$ being connected and $\partial D\not=\emptyset$, given an $f$ on $\partial D$, there exists a continuous function $\widetilde u(x)$ on $D$, such that $\Delta_{((X\cap O)\cup\{x\},\widetilde\mu)}\widetilde u(y)=0$ for $y\in (X\cap O)\cup\{x\}$, where $\widetilde\mu$ is induced by Lemma \ref{bdharmo}. Denote $\omega^D_x(f):=\widetilde u(x)$. Then the topological space $((X,\mu),\tau)$ is a Hausdorff space, where $\tau$ is induced by the topology on the compact surface the graph $(X,\mu)$ embedded on. By Theorem \ref{EU}, Lemma \ref{Addp} and Lemma \ref{bdharmo} we get the following theorem.

\vspace{0.2cm}

\begin{thm}
(1) $\{\omega_\cdot^D\}$ is a harmonic measure system on $(X,\mu)$;
(2) Any open set in $(X,\mu)$ is regular;
(3) The harmonic measure system $\{\omega_\cdot^D\}$ is regular closed;
(4) Strong maximum principle holds for any connected open set with nonempty boundary in $(X,\mu)$.
\end{thm}

\vspace{0.2cm}

\section{Dahlberg's theory}\label{Sec4}

In this section, we generalize Dahlberg's theory \cite{Dahl79} under the axiomatized harmonic measure theory. The proof is based on \cite{Hunt70, Dahl79}.

\vspace{0.2cm}

\subsection{Main theorem}

\vspace{0.2cm}

Suppose $\{\omega_\cdot^{D_\lambda}\}$ is a Radon harmonic measure system on $\mathbb{R}^n$. In this section the Radon harmonic measure system $\{\omega_\cdot^{D_\lambda}\}$ satisfies the following assumptions (H):

\vspace{0.2cm}

\begin{enumerate}
\item[(H1)] $\{D_\lambda\}=\{\text{All the bounded open set in }\mathbb{R}^n\}$.

\vspace{0.2cm}

\item[(H2)] Any bounded Lipschitz domain $D\subseteq\mathbb{R}^n$ is regular relative to $\omega_\cdot^D$.

\vspace{0.2cm}

\item[(H3)] The Radon harmonic measure system $\{\omega_\cdot^{D_\lambda}\}$ is translation and scaling closed \& invariant (TSCI), see the assumption TSCI in Section \ref{sec2.4}.

\vspace{0.2cm}

\item[(H4)] For any bounded open $D\subseteq\mathbb{R}^n$, $\frac{\mathrm{d}\omega^D_x}{\mathrm{d}\omega^D_y}$ is continuous for any $x,y\in D$.
\end{enumerate}

\vspace{0.2cm}

\begin{rem}
(a) Under the assumptions (H2) and (H3), by Theorem \ref{R_Comp}, $\mathbb{R}^n$ is regular closed. By Theorem \ref{Str_MP_G} strong maximal principle is available. (b) By the discussion of Section \ref{sec2.4}, under the assumption (H4) Harnack inequality is available.
\end{rem}

\vspace{0.2cm}

A Lipschitz domain $D$ is called \emph{starlike} about $P_0$ if each ray emanating from $P_0$ intersects $\partial D$ exactly once and the local coordinate system can be taken with its $y$-axis along $QP_0$ for any $Q\in\partial D$.

\vspace{0.2cm}

For any $P\in\partial D$, there exists an open cycle cylinder centered at $P$ with small radius $r>0$,  height $sr>0$ and axis along line $P_0 P$, denoted by $C_{P,r} = C_{P,r,sr}$. The parameter $s$ depends on the the Lipschitz constant $L$ of $\partial D$ such that cylinder $C_{P,r}$ satisfies one bottom is contained in $D$ and the other is outside. For example, we can take $s=2L$.  Denote $\Delta_{P,r}:=C_{P,r}\cap\partial D$. The center point of the bottom of $C_{P,r}$ inside $D$ is denoted by $T_{P,r}$.

\vspace{0.2cm}

The main result of this section is the following theorem.

\vspace{0.2cm}

\begin{thm}\label{thm3.1}
Suppose the Radon harmonic measure system $\{\omega_\cdot^{D_\lambda}\}$ satisfies (H), $D\subset\mathbb{R}^n\,(n\geq 3)$ is a bounded Lipschitz domain, starlike about the origin $O\in D$. If $\mu$ is a positive measure on $D$ and $p >1$, then the following conditions are equivalent:
\begin{enumerate}
\item[(i)] $\mu$ is a Carleson measure in $D$, that is, $\exists\,M>0, \forall\,P\in \partial D, \forall\, r>0$ such that
\begin{equation*}
\mu(C_{P,r}\cap D)\leq M\omega_O^D(\Delta_{P,r});
\end{equation*}
\item[(ii)] $H:f\mapsto\omega_\cdot^D(f)$ is strongly $(p,p)$, that is, $\exists\,K>0,\forall\,f\in L^p(\omega_O^D)$ such that
\begin{equation*}
\int_D|(Hf)(x)|^p\,\mathrm{d}\mu(x) \leq K\int_{\partial D}|f|^p\,\mathrm{d}\omega_O^D.
\end{equation*}
\end{enumerate}
\end{thm}

\vspace{0.2cm}

\subsection{The proof of Theorem \ref{thm3.1}: (ii) $\Rightarrow$ (i)}





\vspace{0.2cm}

\begin{lem}\label{lem3.1a}
	Suppose $D$ is a Lipschitz domain. Then there exists a positive constant $\delta=\delta(s,D)\leq\frac{1}{2}$ such that
	\begin{equation*}
		\omega_{x}^D(\Delta_{P,r})\geq\delta
	\end{equation*}
	for any $x\in C_{P,\delta r}\cap D$.
\end{lem}

\vspace{0.2cm}

\begin{proof}
	Let $\Delta'$ be the bottom of $C_{P,r}$ outside $D$. Then $\chi_{\Delta_{P,r}}(x)\geq\omega_x^{C_{P,r}}(\Delta')$ for all $x\in\partial(D\cap C_{P,r})$. By the strong maximum principle \ref{Str_MP_G}, we get 
	\begin{equation*}
		\omega_x^D(\Delta_{P,r})\geq\omega_x^{C_{P,r}}(\Delta'),
	\end{equation*} 
	for any $x\in C_{P,r}\cap D$. Take $\delta_0:=\inf\limits_{x\in C_{P,\frac{1}{2}r}\cap D}\omega_x^{C_{P,r}}(\Delta')$. By the strong maximum principle \ref{Str_MP_G}, it is clear that $\delta_0<1$, and $\delta_0$ depends only on $s$ and $D$ by the translation and scaling invariant property of the harmonic measure system $\{\omega_\cdot^D\}$ and for any $x\in C_{P,\frac{1}{2}r}\cap D$,
	\begin{equation*}
		\omega_x^D(\Delta_{P,r})\geq\omega_x^{C_{P,r}}(\Delta')\geq\delta_0.
	\end{equation*} 
	Take $\delta:=\min\{\frac{1}{2},\delta_0\}$. Then for any $x\in C_{P,\delta r}\cap D$,
	\begin{equation*}
		\omega_x^D(\Delta_{P,r})\geq\delta,
	\end{equation*}
	where $\delta$ depends only on $s$ and $D$.
\end{proof}

\vspace{0.2cm}

\begin{proof}[The Proof of (ii) $\Rightarrow$ (i)] For any $p>1, P\in \partial D, r>0$ and $x\in C_{P,r}\cap D$, $\omega^D_x(C_{P,\frac{r}{\delta}}\cap\partial D)\geq\delta$, where $\delta\leq\frac{1}{2}$ is the one in Lemma \ref{lem3.1a}. Then by (ii) we get
	\begin{equation*}
		\begin{aligned}
			\delta^p \mu(C_{P,r}\cap D)&\leq\int_{C_{P,r}\cap D}\left(\omega_x^D(C_{P,\frac{r}{\delta}}\cap\partial D)\right)^p\,\mathrm{d}\mu(x)\\
			&=\int_{C_{P,r}\cap D}\left(\int_{\partial D}\chi_{C_{P,\frac{r}{\delta}}\cap\partial D}\,\mathrm{d}\omega_x^D\right)^p\,\mathrm{d}\mu(x)\\
			&\leq \int_{D}\left(\int_{\partial D}\chi_{C_{P,r}\cap\partial D}\,\mathrm{d}\omega_x^D\right)^p\,\mathrm{d}\mu(x)\\
			&\leq K\int_{\partial D}\left(\chi_{C_{P,r}\cap\partial D}\right)^p\,\mathrm{d}\omega^D_O\\
			&=K\omega_O^D(C_{P,r}\cap\partial D).
		\end{aligned}
	\end{equation*}
	Hence $\mu$ is a Carleson measure with its norm $M\leq\delta^{-p}K$.
\end{proof}

\vspace{0.2cm}

\subsection{The proof of Theorem \ref{thm3.1}: (i) $\Rightarrow$ (ii)}

\vspace{0.2cm}

Let $D$ be a bounded Lipschitz domain in $\R^n$.

\vspace{0.2cm}

\begin{prop}[Harnack principle in Lipschitz domain (the inner cone case)]\label{Harnack_Lip}
Let $\Gamma$ be an inner cone with its vertex at $Q\in \partial D$ and $\widetilde D$ be a compact subset in $D\cap \Gamma$. Suppose the harmonic measure system satisfies (H). Then there is a constant $C=C(\widetilde D,\Gamma)$ such that $H(\nu(\widetilde D),D) \leq C$ for any $k>0$ and $Q'\in\partial D$ provided $\nu(D\cap\Gamma)\subset D$, where $\nu(x):=Q'+k(x-Q)$.
\end{prop}

\vspace{0.2cm}

\begin{lem}\label{lem3.2}
If $r$ and $r_0$ with $r\geq r_0>0$ are sufficiently small, depending on the geometry of D, then for any $Q\in\partial D$ and any $P\in D\backslash C_{Q,r}$, 
\begin{equation*}
\omega^D_P(\Delta_{Q,r_0}) \leq C \omega^D_{T_{Q,r}}(\Delta_{Q,r_0}),
\end{equation*}
where $C$ is a constant depending only on $D$.   
\end{lem}

\vspace{0.2cm}

\begin{proof}
Denote $\Delta_r := \Delta_{Q,r} = C_{Q,r} \cap \partial D$ and $T_r := T_{Q,r}$. Then by Lemma \ref{lem3.1a} we get
\begin{equation*} 
\omega^D_{T_{\delta r_0}} (\Delta_{r_0}) \geq \delta,
\end{equation*}
for $0<\delta\leq\frac{1}{2}$, where $\delta$ only depends on $D$. By Harnack principle, 
\begin{equation*}
\omega^D_{T_{r_0}} (\Delta_{r_0}) \geq c_0,
\end{equation*}
where $c_0$ depends on $\delta$ and $D$, hence only depends on $D$. We can choose a constant $c_1 = \frac{1}{c_0}$ such that
\begin{equation*}
\omega_P^D(\Delta_{r_0}) \leq 1\leq c_1 \omega_{T_{r_0}}^D (\Delta_{r_0})
\end{equation*}
for any $P\in D\backslash C_{r_0}$. Suppose that for some positive integer $j$ with $2^j r_0\leq r$,
\begin{equation*}
\omega_P^D(\Delta_{r_0}) \leq c_1 \omega_{T_{2^{j-1}r_0}}^D (\Delta_{r_0})\
\end{equation*}
for any $P\in D\backslash C_{2^{j-1}r_0}$.
	
Let $P_2:=(X_2,0)$, $\Omega =  \{ (X,y) : |X|>1, y>-B|X-X_2|, \sqrt{X^2+y^2} < \frac{diam(D)}{r_0} \}$ where $X_2 = (2,0,\cdots,0)$ in the standard coordinate system. Let $k(P)$ be the harmonic measure in $\Omega$ on the cylinder $\{(X,y)\in\partial\Omega:|X|=1\}$. For fixed point $Q'$ on $\partial C_{2^j r_0} \cap \partial D$, scale $\Omega$ by the factor $2^{j-1}r_0$, then translate it such that the $y$-axis coincides with $OQ$ and $P_2$ becomes $Q'$. We assume that the constant $B$ has been chosen depend on $D$ so that the cone part of $\partial \Omega$ is contained in $\mathbb{R}^n\backslash D$. Let $\Omega_j$ denote the region $\Omega$ scaled by the factor $2^{j-1}r_0$ and let $k_j$ be the corresponding harmonic measure $k$. By the maximum principle,
\begin{equation*}
\omega_P^D(\Delta_{r_0}) \leq k_j(P)
\end{equation*}
for any $P\in D\cap \Omega_j$. And by the strong maximum principle
\begin{equation*}
\omega_P^D(\Delta_{r_0}) \leq  c_1 \omega_{T_{2^{j-1} r_0}}^D(\Delta_{r_0})  k_j(P)
\end{equation*}
for any $P\in \partial (D\cap \Omega_j)$, hence for any $P\in D\cap \Omega_j$ by the maximum principle. 
	
By Harnack principle, there is a constant $c_2$ independent of $j$ such that
\begin{equation*}
\omega_{T_{2^{j-1}r_0}}^D(\Delta_{r_0}) \leq c_2 \omega_{T_{2^{j}r_0}}^D(\Delta_{r_0}).
\end{equation*}
Since $k(P) \to 0$ as $\Omega\ni P\to P_2$, there is an $\eta>0$ such that $k(P) \leq \frac{1}{c_2}$ if $|P-P_2|<\eta$. Hence
\begin{equation}\label{P11}
\omega_P^D(\Delta_{r_0}) \leq  c_1 \omega_{T_{2^{j}r_0}}^D(\Delta_{r_0})
\end{equation}
for any $P\in D\cap\Omega_j$ with $|P-Q'|\leq \eta 2^{j-1}r_0$, hence for any $P\in D\cap \partial C_{2^{j}r_0}$ with $\text{dist}(P,\partial D\cap \partial C_{2^j r_0}) \leq \eta 2^{j-1}r_0$. At the same time, by Harnack inequality there is a constant $c_3$ such that \begin{equation}\label{P12}
\omega_P^D(\Delta_{r_0}) \leq  c_3 \omega_{T_{2^{j}r_0}}^D (\Delta_{r_0})
\end{equation}
for $P\in D\cap \partial C_{2^{j}r_0}$ with $\text{dist}(P,\partial D\cap \partial C_{2^j r_0}) > \eta 2^{j-1}r_0$. Take $c_4 = \max\{c_1,c_3\}$. By \eqref{P11} and \eqref{P12} we get
\begin{equation*}
\omega_P^D(\Delta_{r_0}) \leq  c_4 \omega_{T_{2^{j}r_0}}^D(\Delta_{r_0})
\end{equation*}
for any $P\in D\cap \partial C_{2^{j}r_0}$, hence for any $P\in D\backslash C_{2^{j}r_0}$ by the maximum principle. 
	
Since $r\geq r_0$ there exists a positive integer $j_0$ such that $2^{j_0} r_0 \geq r \geq 2^{j_0-1} r_0$ the conclusion hods by the above inequality and Harnack principle.	
\end{proof}

\vspace{0.2cm}

\begin{lem}\label{lem3.3}
Let $E\subseteq \Delta := B_{Q,r}\cap\partial D$ and $A:=A_{Q,r}$ be the point on the segment $OQ$ whose distance from $Q$ is $r$. Then
\begin{equation*}
\omega^D_A(E) \leq c \frac{\omega^D_O(E)}{\omega^D_O(\Delta)},
\end{equation*}
where $c>0$ depends only on $D$.
\end{lem}

\vspace{0.2cm}

\begin{proof}
	Consider a cone $\Gamma_1$ with vertex at the origin $O$, axis along $OQ$ and the disc $\Delta$ is contained in $\Gamma_1$. Let $H:=H(r,Q) := \Gamma_1 \cap B_{Q,hr}$, where $h$ is a fixed constant. Let $\alpha$ denote the part of $\partial H$ which is contained in the standard inner cone $\Gamma_2$ at $Q$, and let $\beta$ denote the rest of $\partial H$ in $D$. For $h$ sufficiently large, depend on the Lipschitz constant of $D$,  we may assume that $O\not\in B_{Q,hr}$ and there is a $t = t(h), 0<t<1$, such that $tP\in \alpha$ for each $P\in\beta$. Let $\tilde{D} = D \backslash H $. By the maximum principle,
\begin{equation}\label{111}
\omega_P^D(\Delta) \leq \omega_P^{\tilde{D}} (\alpha) + \omega_P^{\tilde{D}} (\beta) \leq (1+t) \omega_P^{\tilde{D}} (\alpha)
\end{equation}
for any $P\in \tilde{D}$. 

Harnack principle implies
\begin{equation*}
\omega_P^D(E) \geq c\omega_A^D(E)
\end{equation*}
for any $P\in\alpha$, where $c$ is independent of $r$. Hence, by the maximum principle,
\begin{equation}\label{112} 
\omega_P^D(E) \geq c\omega_A^D(E) \omega_P^{\tilde{D}} (\alpha)
\end{equation}
for any $P\in \tilde{D}$. 

Setting $P=O$ in \eqref{111} and \eqref{112}, we obtain
\begin{equation*}
\omega_O^D(E) \geq c\omega_A^D(E) \omega_O^{\tilde{D}} (\alpha) \geq \frac{c}{1+t}\omega_A^D(E) \omega_O^D(\Delta)\geq \frac{c}{2}\omega_A^D(E) \omega_O^D(\Delta),
\end{equation*}
which proves the conclusion.
\end{proof}

\vspace{0.2cm}

For any domain $D$, a cone $\Gamma$ with vertex $Q\in\partial D$ is called a nontangential cone at $Q$ if there exists a cone $\widetilde\Gamma$ such that $\bar\Gamma\backslash\{Q\}\subset\widetilde\Gamma\subset D$. 
Denote $k(P,Q) := k^D_O(P,Q) := \frac{\d \omega_P^D(Q)}{\d \omega_O^D(Q)}$, for all $P\in D, Q\in\partial D$.

\begin{lem}\label{lem3.4}
Let $\Gamma$ be any nontangential cone in $D$ at $Q_0$. For any $A'\in\Gamma\cap D$ with $\mathrm{dist}(A',Q_0) = a$, denote $\Delta_j := \Delta_{Q_0,2^j a}$, $j= 0,\cdots,N$, $R_0:=\Delta_0$ and $R_j := \Delta_j\backslash\Delta_{j-1},\, j=1,\cdots,N$. Then
\begin{equation*}
\mathop{\mathrm{ess} \sup}\limits_{Q\in R_j} k(A', Q) \leq \frac{c c_j}{\omega^D_{O}(\Delta_j)},\, j=0,\cdots,N,
\end{equation*}
with $\sum_{j=0}^{N} c_j \leq c' <\infty$, where $c'$ depends on $D$ and $c$ depends on $D$ and the shape of $\Gamma$.  
\end{lem}

\vspace{0.2cm}

\begin{proof}
Let $A$ be the point on the segment $OQ_0$ whose distance from $Q_0$ is $a$. By Harnack priciple $\omega^D_{A'}(E) \leq c \omega^D_A(E)$ where $c$ depends only on $D$ and the shape of $\Gamma$. Hence it is enough to prove the lemma if we show 
\begin{equation}\label{Tra_A}
\frac{\omega^D_A(E)}{\omega^D_O(E)} \leq \frac{c c_j}{\omega^D_{O}(\Delta_j)}
\end{equation} 
for arbitrarily small discs $E \subset R_j $.
	
Let $A_j$ be the point on the segment $OQ_0$ whose distance from $Q_0$ is $2^j a$, then by Lemma \ref{lem3.3} we get 
\begin{equation*}
\omega_{A_j}(E) \leq c\frac{\omega^D_O(E)}{\omega^D_O(\Delta_j)}
\end{equation*}
for any $E \subset \Delta_j$, where $c$ depends only on $D$. Suppose $E$ is small disc with center $Q_j\in R_j, j\geq 4$. By Lemma \ref{lem3.2}, we get
\begin{equation*} 
\omega^D_P(E) \leq c\omega^D_{T_{Q_j,2^{j-2}a}}(E)
\end{equation*}
for any $P\in D\backslash C_{Q_j,2^{j-2}a}$. Harnack principle implies 
\begin{equation*}
\omega^D_{T_{Q_j,2^{j-2}a}}(E) \leq c \omega^D_{A_j}(E).
\end{equation*}
Hence
\begin{equation*} 
\omega^D_P(E) \leq c \frac{\omega^D_O(E)}{\omega^D_O(\Delta_j)}
\end{equation*}
for any $P\in D\backslash C_{Q_j,2^{j-2}a}$. 
	
Let 
\begin{equation*}
\Sigma := \{ (X,y)\in \R^n : |X|^2+y^2 < 1, y>-B|X| \}
\end{equation*}
and $h(P)$ be the harmonic measure of $\{(X,y)\in\partial\Sigma: |X|^2+y^2 =1\}$. Denote $\Sigma_j$ be the region obtained by scaling $\Sigma$ by the factor $2^{j-2}a$ and translating it such that the $y$-axis coincides with $OQ_0$ and the origin coincides with $Q_0$. Denote $h_j$ be the harmonic measure corresponding to $h$. We assume $B$ has been chosen so that the standard outer cone at $Q_0$ coincide with the cone $\{y<-B|X|\}$. Since $D\cap\Sigma_j \subset D\backslash C_{Q_j,2^{j-2}a}$, and the maximum principle imply 
\begin{equation*}
\omega^D_P(E) \leq c \frac{\omega^D_O(E)  }{\omega^D_O(\Delta_j)}h_j(P)
\end{equation*}
for $P\in D\cap \Sigma_j$. Taking $P=A$ we get
\begin{equation*}
\frac{\omega^D_A(E)}{\omega^D_O(E)} \leq \frac{c h_j(A)}{\omega^D_{O}(\Delta_j)}.
\end{equation*} 
For $j=0,1,2,3$, we can choose $c_j$  depend on $D$ by Harnack principle; For $j= 4,\cdots, N$, we can choose $c_j := h_j(A)$. Let $r = \max \{ h(P) : P\in \Sigma, |P| = \frac{1}{2} \}$, then $0<r<1$. by the maximum principle, 
\begin{equation*}
h(\frac{P}{2}) \leq r h(P)
\end{equation*}
for any $ P\in \Sigma$. Hence for $j\geq4$,
\begin{equation*}
h_j(A)=h((0,\frac{1}{2^{j-2}})) \leq r h((0,\frac{1}{2^{j-3}}))= r h_{j-1}(A)\leq\cdots\leq r^{j-4}h_4(A).
\end{equation*}
Since $0<r<1$ we get
\begin{equation*}
\sum_{j=4}^{N} c_j =\sum_{j=4}^{N} h_j(A) \leq \sum_{j=4}^{N} r^{j-4}h_4(A)\leq \frac{1}{1-r}<\infty.
\end{equation*} 
Then \eqref{Tra_A} is proved, which implies the conclusion of the lemma.	
\end{proof}

\vspace{0.2cm}

\begin{lem}\label{lem3.5}
Denote $\Delta := \Delta_{Q_0,r_0}$ for $Q_0\in \partial D$. Then
\begin{equation*}
\mathop{\mathrm{ess} \sup}\limits_{Q'\in \partial D \backslash \Delta} k(P,Q')\to 0,
\end{equation*}
as $P\to Q_0$. 
\end{lem}

\vspace{0.2cm}

\begin{proof}
Denote $\Sigma_{0}$ be the region by Scaling $\Sigma$ in the proof of lemma \ref{lem3.4} by the factor $\frac{r_0}{3}$ and translating it to $Q_0$ as in the proof of lemma \ref{lem3.4}. Denote $h_0$ be the harmonic measure corresponding to $h$ defined in the proof of Lemma \ref{lem3.4}. For any small disc $\Delta' := \Delta'_{Q',r'}$ with $Q'\in \partial D\backslash\Delta$ and $0<r'<\frac{r_0}{3}$, Lemma \ref{lem3.2} implies 
\begin{equation*}
\omega^D_P(\Delta') \leq c_1 \omega^D_{T_{Q',\frac{r_0}{3}}}(\Delta')
\end{equation*}
for any $P\in D\backslash C_{Q',\frac{r_0}{3}}$, where $c_1$ depends only on $D$. By Harnack inequality, we get
\begin{equation*}
\omega^D_{T_{Q',\frac{r_0}{3}}}(\Delta')\leq c_2\omega^D_O(\Delta'),
\end{equation*}
where $c_2$ depends on $D$ and $r_0$. Then the above two inequalities imply
\begin{equation*}
\omega^D_P(\Delta') \leq c_1c_2 \omega^D_O(\Delta')
\end{equation*}
for any $P\in \partial\Sigma_0\cap D\subseteq D\backslash C_{Q',\frac{r_0}{3}}$. Since $0\leq \omega_P(\Delta') \leq 1 = h_0(P)$ for $P\in\partial \Sigma_0 \cap D$ and $ \omega_P(\Delta') =0\leq h_0(P)$ for $P\in\partial D\cap \Sigma_0$, the maximum principle implies
\begin{equation*}
\omega^D_P(\Delta') \leq h_0(P),
\end{equation*}
for any $P\in \Sigma_0\cap D$. Hence 
\begin{equation*}
\omega^D_P(\Delta') \leq c_1c_2 \omega^D_O(\Delta') h_0(P)
\end{equation*}
for any $P\in\partial (\Sigma_0\cap D)$, then for $P\in\Sigma_0\cap D$ by the maximum principle. Since $h_0(P)\to 0$ as $P$ tends to $Q_0$, we have
\begin{equation*}
\frac{\omega^D_{P}(\Delta')}{\omega^D_O(\Delta')} \to 0
\end{equation*}
as $P\to Q_0$ uniformly for any radius small $\Delta'\subset\partial D\backslash D$ containing $Q'$. This proves the conclusion by the definition of $k(P,Q')$.
\end{proof}

\vspace{0.2cm}

We need to prove the following lemma.

\vspace{0.2cm}

\begin{lem}\label{lem3.6}
Suppose $D\subset \R^n (n\geq 3)$ is a Lipschitz domain, starlike about the origin $O\in D$ and $\mu$  is a Carleson measure in $D$. For a $Q_0\in \partial D$ and appropriate $A,B>0$ depend on $D$, denote $S_1 := \{Q\in \partial D: \mathrm{dist}(Q,l_{OQ_0})<B\}$, $D_2 := \{P\in D: \mathrm{dist}(P,Q_0)<2A, \mathrm{dist}(P,l_{OQ_0})<2B)\}$. Then for any $f\in L^p(\partial D,\omega_\cdot^D)$ and any $p\geq 1$, there exists $K>0$ such that
\begin{equation*}
\mu(\{ P\in D_2: |\omega^D_P(f)|> s \}) \leq K s^{-p}\|f\|^p_{p,\omega^D_O}.
\end{equation*}
\end{lem}

\vspace{0.2cm}

\begin{proof}
It is sufficient to prove the lemma in the case $f\geq 0$, since $\mu\{ P\in D_2: |\omega^D_P(f)|> s \} \leq \mu\{ P\in D_2: |\omega^D_P(f^+)|> \frac{s}{2} \} + \mu\{ P\in D_2: |\omega^D_P(f^-) |> \frac{s}{2} \} $.
	
Denote $S_i:= \{Q\in \partial D: \mathrm{dist}(Q,l_{OQ_0})<iB\}$ and $D_i := \{P\in D: \mathrm{dist}(P,Q_0)<iA, \mathrm{dist}(P,l_{OQ_0})<iB\}$ for $i\in\mathbb{Z}^+$. For any $Q\in S_3$ we can choose $\Gamma(Q)$ being a small nontangential inner cone with center line paralleling to $OQ_0$ and the height $h>0$ small enough such that: 
\begin{enumerate}
\item[(1)] For any $Q\in S_2$, $\Gamma(Q)\subseteq D_3$ and
\item[(2)] by Lemma \ref{lem3.5}, we can choose $\Delta_h :=\Delta_{Q,r_h}$ such that for any $Q\in S_3$ and any $P'\in\Gamma(Q)$,
\begin{equation*}
\mathop{\mathrm{ess} \sup}\limits_{Q'\in\partial D\backslash \Delta_h}k(P',Q') \leq cc',
\end{equation*}
where $c, c'$ are the constants in Lemma \ref{lem3.4}.
\end{enumerate} 
Then Lemma \ref{lem3.4} implies that for any $Q\in S_3$ and any $P'\in \Gamma(Q)$,
\begin{equation*}
\begin{aligned}
\omega^D_{P'}(f) 
&= \int_{\Delta_h} f\,\mathrm{d}\omega^D_{P'} +\int_{\partial D\backslash\Delta_h} f\d \omega^D_{P'}\\
&\leq \sum_{j=0}^{m(|P'Q|)}\int_{R_j} fk(P',\cdot)\d \omega^D_{O} + \int_{\partial D\backslash \Delta_h} fk(P',\cdot)\d \omega^D_{O} \\
&\leq c {\sum_{j=0}^\infty {c_j}\fint_{\Delta_j} f\d \omega^D_{O}} + \mathop{\mathrm{ess} \sup}\limits_{Q'\in\partial D\backslash \Delta_h}k(P',Q') \fint_{\partial D} f \d \omega^D_{O} \\ 
&\leq c \sum_{j=0}^\infty {c_j}f^*(Q) + cc'f^*(Q) \\ 
&= 2cc'f^*(Q), 
\end{aligned}
\end{equation*}
where $c'$ depends only on $D$, $c$ depends on $D$ and the shape of $\Gamma(Q)$, hence depends only on $D$ and 
\begin{equation*}
f^*(Q) := \sup\limits_{r>0} \fint_{\Delta(Q,r)}f \d\omega^D_{O}.
\end{equation*} 
Then we get for any $Q\in S_3$,
\begin{equation}\label{equ3.1}
\sup_{P'\in \Gamma(Q)}\omega^D_{P'}(f)\leq C_0 f^*(Q), 
\end{equation}
where $C_0$ depends only on $D$.

\vspace{0.2cm}

We divide the rest of proof into two steps.

\vspace{0.2cm}
	
Step 1. Let $V := D_2 \cap (\mathop{\cup}\limits_{Q\in S_2} \Gamma(Q))$, $U := D_2-V$. By the Harnack inequality and H\"older inequality, $\omega^D_P(f) \leq C \omega^D_{O}(f)\leq C \|f\|_{p,\omega_O^D}$ for any $P\in U$ and $f\in L^p(\partial D,\omega_\cdot^D)$. Hence, for $s\geq C^{\frac{1}{p}}\|f\|_{p,\omega_O^D}$, we get $C \omega^D_{O}(\left(\frac{f}{s}\right)^p)\leq 1$, which implies $\omega^D_P(\frac{f}{s})\leq 1$, i.e.  
\begin{equation*}
\mu \{P\in U: \omega^D_P(f) >s \}=\mu \{P\in U: \omega^D_P(\frac{f}{s}) >1 \}=0\leq C s^{-p} \|f\|_{p,\omega^D_O}^p;
\end{equation*}
And for $0<s< C^{\frac{1}{p}}\|f\|_{p,\omega_O^D}$, which implies $1<Cs^{-p} \|f\|_{p,\omega^D_O}^p$, then
\begin{equation*}
\mu \{P\in U: \omega^D_P(f) >s \}\leq \mu(D)\cdot1< \mu(D)Cs^{-p} \|f\|_{p,\omega^D_O}^p.
\end{equation*}
Then the above two inequalities imply for any $s>0$,
\begin{equation}\label{equ3.2}
\mu \{P\in U: \omega^D_P(f) >s \}\leq C' s^{-p} \|f\|_{p,\omega^D_O}^p,
\end{equation}	
where $C':=\max\{C,\mu(D)C\}$ depends only on $D$.

\vspace{0.2cm}

Step 2. We will prove 
\begin{equation}\label{4.888}
	\mu \{ P\in V: \omega_P^D(f) >s \} \leq C s^{-p} \|f\|_{p,\omega_O^D}^p.
\end{equation}

Let $P^*\in \partial D$ be the projection of $P\in D$ along the direction of $OQ_0$. Define  $R(P) := \{ P'\in \R^n: \mathrm{dist}(P',P) < 2\mathrm{dist}(P,P^*), \mathrm{dist}(P',l_{PP^*}) < t\, \mathrm{dist}(P,P^*) \}$, where $t>0$ is chosen so small that, for any $P\in V$, $P \in \Gamma(Q)$ for any $Q\in R(P)\cap S_3$. Moreover, if $\omega^D_P(f) > s$ for $P\in V$, then $f^*(Q)>\frac{1}{C_0}s$ for any $Q\in R(P)\cap S_3$ by estimate \eqref{equ3.1}. 

Suppose $F\subset \{P\in V: \omega_P^D(f) >s \}$ is compact. Let $R'(P)$ be the maximum rectangle in $R(P)$. The radius of $R'(P)$ is denoted by $r(P): = \frac{\sqrt{2}}{2}t\mathrm{dist}(P,P^*)$. Let $r_0:=\mathop{\sup}\limits_{P\in F}r(P)$. Then $r_0\leq \frac{\sqrt{2}}{2}th$. Take $P_1\in F$ such that $r(P_1)>\frac{r_0}{2}$ and define
\begin{equation*}
	F_1:= F\backslash R'(P_1),\, r_1:= \sup_{P\in F_1} \{ r(P)\}.
\end{equation*}
Take $P_2\in F_1$ such that $r(P_2) >\frac{r_1}{2}$ and define
\begin{equation*}
	F_2:= F\backslash (R'(P_1)\cup R'(P_2)),\, r_2:= \sup_{P\in F_2} \{ r(P)\}.
\end{equation*}
Continue this process until the first integer $m$ is found such that $F_m$ is an empty set. If such an integer does exists, continue this process and set $m=+\infty$. If this case occurs,  then $F\backslash\mathop{\cup}\limits_{i=1}^\infty R'(P_i)\not=\emptyset$.

We claim that $\{R'(P_j)\}$ we constructed as the above is a finite cover of $F$. Suppose $i>j$. Then $P_i\in F_{i-1}\subset F\backslash R'(P_j)$, thus $P_i\notin R'(P_j)$ and $r(P_i)\leq r_{i-1} <r_{j-1} <2r(P_j)$. It easy to get $\frac{1}{3}R'(P_i) \cap \frac{1}{3}R'(P_j) = \emptyset$ by $P_i\notin R'(P_j)$, where $\frac{1}{3}R'(P) := R'(\frac{2}{3}P^*+\frac{1}{3}P)$. If $m=+\infty$, then there exists a subsequence of $\{r(P_i)\}$ converging to zero since $\frac{1}{3}R'(P_i)$ are pairwise disjoint and have centers in a bounded set $\frac{1}{3}F\subset V$. In this case, since $F\backslash\mathop{\cup}\limits_{i=1}^\infty R'(P_i)\not=\emptyset$, there exists a $P$ satisfying $\frac{1}{3}P \in \frac{1}{3}F\backslash\cup_{i=1}^\infty \frac{1}{3}R'(P_i)$. Then $P\in F_i$ and $r(P)<r(P_i)$ for all $i\in \N$. Hence $r(P)=0$ and $R(P)=\emptyset$, which is a contradiction. Thus $m < \infty$.

We claim that $\sum\limits_{j=1}^m \chi_{R'(P_j)} \leq 2^n$.  Fix a $y\in \R^n$, then $\R^n$ can be decomposed into $2^n$ parts by $n$ hyperplanes which are parallel to coordinate hyperplanes and pass through the point $y$. Choose one of these $2^n$ parts and denote it by $O$. Suppose $P_i, P_j \in O$ $(1\leq i,j\leq m)$ such that $y\in R'(P_i)\cap R'(P_j)$ and $r(P_i)\geq r(P_j)$. Then $P_j\in R'(P_i)$, contradicting the fact that each $P_j$ belongs to no other $R'(P_i)$. Hence there exists at most one rectangle $R'(P_i)$ with center in $O$ containing $y$. Then the conclusion of the claim holds.

By the above two claims, since $\mu$ is a Carleson measure and $\omega_{P_j}^D(f)>s$ which implies $f^*(Q)>\frac{1}{C_0}s$ for any $Q\in R(P_j)\cap S_3$, we get 
\begin{equation*}
\begin{aligned}
\mu(F)&\leq \sum_{j=1}^m \mu(R'(P_j))\leq \sum_{j=1}^m M  \omega^D_O(R(P_j)\cap \partial D) \\
&\leq  2^n M \omega^D_O(\{ Q\in S_3:f^*(Q)>\frac{1}{C_0} s \}).
\end{aligned}
\end{equation*} 
Since $F$ is arbitrary, it follows that 
\begin{equation}\label{equ3.3} 
\mu(\{P\in V:\omega^D_P(f) >s\})\leq C_1\omega^D_O(\{Q\in S_3:f^*(Q)>\frac{1}{C_0} s\}), 
\end{equation}
where $C_1:=2^n M$.

Moreover, by the H\"older inequality, for $Q\in S_3$, we get 
\begin{equation*}
\fint_{\Delta_{Q,r}} fd\omega_O^D \leq (\fint_{\Delta_{Q,r}}f^{p}d\omega_O^D)^{\frac{1}{p}} (\fint_{\Delta_{Q,r}}1^{p'}d\omega^D_O)^{\frac{1}{p'}}=(\fint_{\Delta_{Q,r}}f^{p}d\omega_O^D)^{\frac{1}{p}}.
\end{equation*}
Hence 
\begin{equation}\label{maf}
f^*(Q) \leq ((f^p)^*(Q))^{\frac{1}{p}}.
\end{equation}
From the ordinary maximal inequality (see \cite{Stei70} p.5) it follows that
\begin{equation}\label{SteinIneq}
\omega_O^D(\{(f^p)^*(Q)> C_0^{-p}s^p\}) \leq \frac{C_2C_0^p}{s^p}\omega_O^D(f^p)= C_2C_0^p s^{-p} \|f\|^p_{p,\omega_O^D},
\end{equation}
where $C_2$ depends only on $n$. By \eqref{equ3.3}, \eqref{maf} and \eqref{SteinIneq} we get
\eqref{4.888}.

\vspace{0.2cm}

Finally, \eqref{equ3.2} and \eqref{4.888} imply the conclusion of the lemma.
\end{proof}

\vspace{0.2cm}

We are now ready to give the proof of (i) $\Rightarrow$ (ii).

\vspace{0.2cm}

\begin{proof}[The proof of (i) $\Rightarrow$ (ii)] 
As in the proof of Lemma \ref{lem3.6}, it is sufficient to treat the case $f\geq 0$. From compactness of $\partial D$ we get the existence of finitely many domains $D_0, D_1, \cdots ,D_N$ with the following properties:
\begin{enumerate} 
\item[(1)] $D_i\subset D$ for $0 \leq i \leq N$ and $D= \mathop{\cup}\limits_{i=0}^N D_i$;
\item[(2)] Each $D_i$ with $1 \leq i \leq N$ is a domain of the type $D_2$ defined in Lemma \ref{lem3.6}. 
\end{enumerate}
Then by Lemma \ref{lem3.6} it follows that for any $q\geq 1$,
\begin{equation*}
\mu\{ P\in D: |\omega^D_P(f)|> s \} \leq C s^{-q} \int_{\partial D }|f|^q d\omega^D_O,
\end{equation*}
which implies $H:f\mapsto\omega_\cdot^D(f)$ is weakly $(q,q)$. At the same time, $H:f\mapsto\omega_\cdot^D(f)$ is weakly $(\infty,\infty)$ by the maximum principle. Then by the Marcinkiewicz interpolation theorem $H:f\mapsto\omega_\cdot^D(f)$ is strongly $(p,p)$ for any $q < p < \infty$. Choose $q=1$, we get the conclusion.
\end{proof}

\section{Acknowledgments}
This research work is partly supported by the National Natural Science Foundation of China (NSFC 12371096). We would like to thank Professor Jun Geng and Sibei Yang for their valuable suggestion for this research work.





\renewcommand{\baselinestretch}{0.1}
\bibliographystyle{plain}

\end{document}